\documentclass{amsart}

\usepackage{amsmath, graphicx, hyperref, cite, enumerate, comment}

\newtheorem{theorem}{Theorem}[section]
\newtheorem{lemma}[theorem]{Lemma}
\newtheorem{proposition}[theorem]{Proposition}

\theoremstyle{definition}
\newtheorem{definition}[theorem]{Definition}

\theoremstyle{remark}
\newtheorem{remark}[theorem]{Remark}

\numberwithin{equation}{section}

\newcommand{\abs}[1]{\left|#1\right|}

\newcommand{\Rm}{\textup{Rm}}
\newcommand{\Ric}{\textup{Ric}}

\newcommand{\FS}{\textup{FS}}
\newcommand{\Tr}{\textup{Tr}}

\newcommand{\eval}[2]{\left. #1 \right|_{#2}}

\newcommand{\R}{\mathbb{R}}
\newcommand{\C}{\mathbb{C}}

\newcommand{\N}{\mathbb{N}}
\newcommand{\Z}{\mathbb{Z}}

\newcommand{\ddt}{\frac{d}{dt}}
\newcommand{\dd}[1]{\frac{\partial}{\partial #1}}
\newcommand{\D}[2]{\frac{\partial #1}{\partial #2}}
\newcommand{\ddbar}{\partial \bar{\partial}}
\newcommand{\barrr}[1]{\overline{#1}}

\begin{document}

\title[KRF on projective bundles over K\"ahler-Einstein manifolds]{K\"ahler-Ricci flow on projective bundles over K\"ahler-Einstein manifolds}

\author[Frederick T.-H. Fong]{Frederick Tsz-Ho Fong}
\address{Department of Mathematics, Stanford University, Stanford, California 94305}
\email{thfong@math.stanford.edu}
\thanks{The author was supported in part by NSF Grant DMS-\#0604960.}


\date{submitted on April 12, 2011, revised on October 11, 2011.}

\begin{abstract}
We study the K\"ahler-Ricci flow on a class of projective bundles $\mathbb{P}(\mathcal{O}_\Sigma \oplus L)$ over compact K\"ahler-Einstein manifold $\Sigma^n$. Assuming the initial K\"ahler metric $\omega_0$ admits a $U(1)$-invariant momentum profile, we give a criterion, characterized by the triple $(\Sigma, L, [\omega_0])$, under which the $\mathbb{P}^1$-fiber collapses along the K\"ahler-Ricci flow and the projective bundle converges to $\Sigma$ in Gromov-Hausdorff sense. Furthermore, the K\"ahler-Ricci flow must have Type I singularity and is of $(\C^n \times \mathbb{P}^1)$-type. This generalizes and extends part of Song-Weinkove's work \cite{SgWk09} on Hirzebruch surfaces.
\end{abstract}

\maketitle

\section{Introduction}
The Ricci flow was introduced by Hamilton in his seminal paper \cite{Ham82} in 1982, proving the existence of constant sectional curvature metric on any closed 3-manifold with positive Ricci curvature. Since then, the Ricci flow has been making breakthroughs in settling several long-standing conjectures. Just to name a few, based on a program proposed by Hamilton, a complete proof of the Poincar\'e conjecture was given by Perelman \cite{Perl1, Perl2, Perl3} around 2003. See also \cite{CaoZh06, KL08, MgTn06}. Furthermore, the Differentiable Sphere Theorem was proved by Brendle-Schoen \cite{BdSch09} in 2007, giving an affirmative answer to a conjecture about differential structures of quarter-pinched manifolds proposed by Berger and Klingenberg in 1960s. In the realm of K\"ahler geometry, the K\"ahler-Ricci flow was introduced by Cao in \cite{Cao85}, which proves the smooth convergence towards the unique K\"ahler-Einstein metric in the cases $c_1 < 0$ and $c_1 = 0$.

There has been much interest in understanding the limit behavior and singularity formation of the Ricci flow in both Riemannian and K\"ahler settings. Hamilton introduced in \cite{Ham95} a method of studying singularity formation of the Ricci flow by considering the Cheeger-Gromov limit of a sequence of rescaled dilated metrics. The singularity model obtained, which is often an ancient or eternal solution, captures the geometry of the singularity formation near the blow-up time of the flow. For closed 3-manifolds, the study of ancient $\kappa$-solutions formed by the dilated sequence limit in Hamilton-Perelman's works (e.g. \cite{Ham95, Perl1}) leads to a solid understanding of singularity formation of closed 3-manifolds. 

Another way of interpreting singularity formation is by the Gromov-Hausdorff limit, regarding the manifold as a metric space. This notion was employed recently in the study of algebraic varieties by Song, Tian, Weinkove et. al in \cite{SgTn06, Tn08, SgTn09, SgWk10, SgWk11, SSW11}. The unified theme of these works is the conjecture that the K\"ahler-Ricci flow will carry out an analytic analogue of Mori's minimal model program which is about searching for birationally equivalent models ``minimal" in some algebraic sense. Like Hamilton-Perelman's work, a surgery may need to be performed in continuing the flow if necessary. To this end, the Gromov-Hausdorff convergence provides a bridge to continue the relevant geometric data.

For a better understanding of singularity formation of the K\"ahler-Ricci flow, one could study some algebraically concrete spaces and explore their flow behavior and possible singularity types and models. In the work by Feldman-Ilmanen-Knopf \cite{FIK03}, Cao \cite{Cao96} and Koiso \cite{Koi}, gradient K\"ahler-Ricci solitons were constructed on the $\mathcal{O}(-k)$-bundles over $\mathbb{P}^n$. Their work employs the $U(n+1) / \mathbb{Z}_k$-symmetry introduced by Calabi in \cite{Cal82} which reduces the K\"ahler-Ricci flow equation to a PDE with one spatial variable. Assuming Calabi's symmetry, Song-Weinkove \cite{SgWk09} characterized the limit behavior (in the Gromov-Hausdorff sense) of the Hirzebruch surfaces $\mathbb{P}(\mathcal{O} \oplus \mathcal{O}(-k))$ and their higher dimensional analogues, which are $\mathbb{P}^1$-bundles over $\mathbb{P}^n$. In their paper, it was proved that the K\"ahler-Ricci flow exhibits three distinct behaviors: (1) collapsing along the $\mathbb{P}^1$-fibers; (2) contracting the exceptional divisor; or (3) shrinking to a point. This trichotomy is determined by the triple $(n, k, [\omega_0])$ where $[\omega_0]$ is the initial K\"ahler class. Later in \cite{SgWk10} by the same authors, case (2) is much generalized and the assumption on the symmetry is removed. The Calabi symmetry assumption is removed in case (1) by a recent preprint \cite{SSW11} by Song, Sz\'ekelyhidi and Weinkove.

The purpose of this paper is two-fold. For one thing, we generalize Song-Weinkove's work \cite{SgWk09} on Hirzebruch surfaces to a class of projective bundles over any compact K\"ahler-Einstein manifold. We will employ an ansatz, known as the momentum construction, which coincides with Calabi's $U(n+1) / \Z_k$-symmetry on Hirzebruch surfaces where the base manifold has the Fubini-Study metric. The idea of the momentum construction of projective bundles was introduced and studied in the subject of extremal K\"ahler metrics by Hwang-Singer in \cite{HgSgr} and by Apostolov-Calderbank-Gauduchon-(T{\o}nnesen-Friedman) in \cite{ACGT}. We will show that under this momentum construction, one can give a cohomological criteria under which the K\"ahler-Ricci flow will collapse the $\mathbb{P}^1$-fiber near the singularity similar to the Hirzebruch surface cases in \cite{SgWk09}. Secondly, we study the singularity model of these projective bundles (including Hirzebruch surfaces) via the techniques developed by Hamilton in \cite{Ham95}. We show that these collapsing projective bundles equipped with momenta will all exhibit $\C^n \times \mathbb{P}^1$-singularities, and also that the Ricci flow solution has a Type I singularity. Here is the summary of our results:\\

\noindent \textbf{Main results.} \label{main} Let $M = \mathbb{P}(\mathcal{O}_\Sigma \oplus L)$ be a projective bundle where $(\Sigma, \omega_\Sigma)$ is a compact K\"ahler-Einstein manifold such that $\Ric(\omega_\Sigma) = \nu\omega_\Sigma$ for some $\nu\in\R$, and $L \to \Sigma$ is a holomorphic line bundle that admits a Hermitian metric $h$ such that the Chern curvature is given by $F_\nabla = -\lambda \omega_\Sigma$, $\lambda > 0$. Let $\omega_0$ be a K\"ahler metric on $M$ constructed by a $U(1)$-invariant momentum profile with K\"ahler class $[\omega_0] = \lambda b_0[\Sigma_\infty] - \lambda a_0 [\Sigma_0]$. Suppose the triple $(\Sigma, L, [\omega_0])$ satisfies the following conditions
\begin{align*}
\nu & \leq \lambda, \text{ or }\\
\nu & > \lambda \text{ and } (\nu - \lambda)b_0 < (\nu + \lambda)a_0,
\end{align*}
then along the K\"ahler-Ricci flow $\partial_t \omega_t = -\Ric(\omega_t)$, $t \in [0,T)$, we have
\begin{itemize}
\item $(M, g(t))$ converges in Gromov-Hausdorff sense to $(\Sigma, \omega_\Sigma)$ (Theorem \ref{GrHdf});
\item the associated ancient $\kappa$-solution is $\C^n \times \mathbb{P}^1$ (Theorem \ref{type_I});
\item the Ricci flow solution must have a Type I singularity (Theorem \ref{no_type_II}).
\end{itemize}

This paper is organized as follows. Sections 2 and 3 are the preliminaries which define our projective bundles and construct K\"ahler metrics using momentum profiles. We will see that the K\"ahler-Ricci flow is equivalent to a heat-type equation for the evolving momentum profile. Section 4 explains the trichotomy of blow-up exhibited by different choice of the triples $(\Sigma, L, [\omega_0])$ via the calculation of K\"ahler classes and Chern classes. Section 5 is a variation on the theme of Song-Weinkove's work \cite{SgWk09} on Hirzebruch surfaces (the collapsing case). We show that similar limiting behavior can be observed in our projective bundles. Sections 6 and 7 are about singularity analysis using rescaled dilations. We show in Section 6 that the ancient $\kappa$-solution obtained from the Cheeger-Gromov limit must split into a product. We will classify their singularity type and the curvature blow-up rate in Section 7.

We also acknowledge that K\"ahler-Ricci solitons on this category of bundles (and their variants) were studied and constructed in \cite{DcrWg11, BYg08, CLi10}.\\

\noindent \textbf{Acknowledgments.} The author would like to express his heartfelt gratitude to his advisor Professor Richard Schoen for all his continuing support and many productive discussions. The author would also like to thank Yanir Rubinstein for arousing his interest in this topic and for many helpful ideas, and also Ziyu Zhang for informing him of some algebraic aspects related to this study.

\section{Projective Bundles}
In this section, we will define and elaborate on the projective bundles under consideration in this paper. We first start with a compact K\"ahler-Einstein manifold $\Sigma^n$ with $\dim_\C = n$. A K\"ahler manifold is called K\"ahler-Einstein if it admits a K\"ahler form $\omega_\Sigma$ whose Ricci form is a real constant multiple of $\omega_\Sigma$, i.e. $\Ric(\omega_\Sigma) = \nu\omega_\Sigma$, $\nu \in \R$. Clearly a necessary condition for a compact K\"ahler manifold to be K\"ahler-Einstein is that the first Chern class $c_1$ has a definite sign. It is well-known by results of Aubin \cite{Au76} and Yau \cite{Yau78} that when $c_1 < 0$ or $=0$ K\"ahler-Einstein metric always exists. However, if $c_1 > 0$ (i.e. Fano manifolds), K\"ahler-Einstein metrics do not exist in general. For compact Riemann surfaces, i.e. $\dim_\C = 1$, K\"ahler-Einstein metric must exist according to the classical uniformization theorem. See also Cheng-Yau's work \cite{ChgYau80} on pseudoconvex domains in the complete non-compact case.

In this article, we will not go into the detail of existence issues of K\"ahler-Einstein metrics, but we will start with a compact K\"ahler manifold $\Sigma^n$ which is equipped with a K\"ahler-Einstein metric $\omega_\Sigma$, such that the Ricci form is given by $\Ric(\omega_\Sigma) = \nu\omega_\Sigma$ where $\nu\in\R$. We take this K\"ahler-Einstein manifold to be our base manifold, and build a projective $\mathbb{P}^1$-bundle upon it. Precisely, we construct our projective bundles as follows:
\[M = \mathbb{P}(\mathcal{O}_\Sigma \oplus L).\]
Here $\mathcal{O}_\Sigma$ is the trivial line bundle, and $L \to \Sigma$ is a holomorphic line bundle which is equipped with a Hermitian-Einstein metric $h$ such that $\sqrt{-1}\ddbar\log h = \lambda \omega_\Sigma$, $\lambda\in\R$. Here $\mathbb{P}$ denotes the projectivization of the holomorphic rank-2 bundle $\mathcal{O}_\Sigma \oplus L$ over $\Sigma$. The local trivialization $(z, u)$ of this rank-2 bundle has transition functions of the form $(z_\alpha, u_\alpha) \approx (z_\beta, \eta_{\alpha \beta} u_\alpha)$ for some $\eta_{\alpha\beta} \in \check{H}^1(\Sigma, \mathcal{O}_\Sigma^*)$. Passing to the projectivization quotient, every element under this trivialization can be expressed as either $[1:\frac{u}{z}]$ for $z \not =0$ or $[0:1]$ and we may regard $[0:1]$ as the infinity. One can check easily that the projectivization factors through the identification by the transition functions $\mathcal{O}_\Sigma \oplus L$. Therefore, one can regard the projectivization of $\mathcal{O}_\Sigma \oplus L$ as compactifying each fiber by adding an infinity point $(x, [0:1])$ and hence $M$ can be regarded as a $\mathbb{P}^1$-bundle over $\Sigma$. We define $\Sigma_0$ to be the zero section $\{x: [1:0]\}$ and $\Sigma_\infty$ to be the infinity section $\{x:[0:1]\}$. It is easy to see that the zero section $\Sigma_0$ and the infinity section $\Sigma_\infty$ are global over $\Sigma$.

The class of holomorphic line bundles over $\Sigma$ with tensor product as the operation form a group which is known as the Picard group, denoted by Pic$(\Sigma)$. For $\Sigma = \mathbb{P}^n$, it is well-known (see e.g. \cite{GH}) that Pic$(\mathbb{P}^n) = \Z$ and the line bundles over $\mathbb{P}^n$ are given by $\mathcal{O}_{\mathbb{P}^n}(k)$, $k \in \Z$. In particular if $(\Sigma, L) = ( \mathbb{P}^1, \mathcal{O}_{\mathbb{P}^1}(-k))$, $k > 0$, the projective bundles $M = \mathbb{P}(\mathcal{O}_{\mathbb{P}^1} \oplus \mathcal{O}_{\mathbb{P}^1}(-k))$ are called the Hirzebruch surfaces. When $k=1$, the projective bundles is $\mathbb{P}^2 \# (-\mathbb{P}^2)$, i.e. $\mathbb{P}^2$ blown-up at a point. When $\Sigma = \C / \Lambda$, i.e. an elliptic curve or a 2-torus, the class of line bundles are classified by a classical result by Appell-Humbert (see \cite{Mfd}). In general for Riemann surface $\Sigma_g$ of genus $g$, the Picard group Pic$(\Sigma_g)$ is isomorphic to $J(\Sigma_g) \times \Z$ where $J(\Sigma_g)$ is a compact complex manifold $\C^g / \Lambda$ of dimension $g$.

The projective bundle $M$ under our consideration is characterized by the pair $(\Sigma, L)$ where $\Sigma$ is a compact K\"ahler-Einstein manifold and $L$ a holomorphic line bundle over $\Sigma$ which is equipped with a Hermitian metric $h$ such that the Chern cuvature is of the form $F_{\nabla} = -\lambda \omega_\Sigma$. In particular, the line bundles generated by $\det(T^*M)$ all fall into this category. Moreover, we will only focus on line bundles $L$ with $\lambda > 0$, since projective bundles $\mathbb{P}(\mathcal{O}_\Sigma \oplus L)$ is biholomorphic to its dual cousin $\mathbb{P}(\mathcal{O}_\Sigma \oplus L^*)$. Since $c_1 (L) = -c_1 (L^*)$, one can replace $L$ by $L^*$ in case $c_1(L)$ is negative. We do not discuss the case of flat bundles, i.e. $\lambda = 0$, in this paper.

\section{U(1)-invariant K\"ahler metrics}
Let's first recapitulate the construction of the category of projective bundles we concerned about in the rest of this article. We let $M = \mathbb{P}(\mathcal{O}_\Sigma \oplus L)$, where $(\Sigma, \omega_\Sigma)$ is a K\"ahler-Einstein manifold such that $\Ric(\omega_\Sigma) = \nu\omega_\Sigma$, $\nu\in\R$. Suppose $L$ is a holomorphic line bundle over $\Sigma$ such that it equips with a Hermitian metric $h$ whose Chern curvature is of the form $F_\nabla = - \lambda\omega_\Sigma$, $\lambda > 0$. In particular, such a Hermitian metric $h$ must exist if $\omega_\Sigma$ is a compact Riemann surface.

We will discuss the construction of $U(1)$-invariant K\"ahler metrics on these projective bundles in this section. Regard the circle group $U(1)$ as $\{e^{i\theta}: \theta \in [0, 2\pi)\}$. The $U(1)$-action defined by
\[e^{i \theta} \cdot (x, [z:u]) = (x, [z: e^{i\theta} u]).\]
Clearly, the action factors through the transition functions of the bundle, and fixes the zero and infinity sections. 

Recall that $\omega_\Sigma$ be the K\"ahler-Einstein form on the manifold $\Sigma$ and we have $\Ric(\omega_\Sigma) = \nu\omega_\Sigma$ for some $\nu \in \R$. Using the Hermitian-Einstein metric $h$ described above, one can define a height parameter $\rho$ on $M \backslash (\Sigma_0 \cup \Sigma_\infty)$ given by
\[\rho = \log \| \cdot \|^2_h.\]
Note that $\rho = -\infty$ corresponds to the zero section and $\rho = \infty$ corresponds to the infinity section.

Our next step is to define K\"ahler metrics on $M$ which is invariant under the circle action defined above. We start by looking for possible K\"ahler classes that $M$ can have. We denote $[\Sigma_0]$ and $[\Sigma_\infty]$ as the Poincar\'e duals (with respect to a fixed background volume form) of $\Sigma_0$ and $\Sigma_\infty$ in $H_2 (M,\R)$ respectively, i.e. $\int_{\Sigma_\infty} [\Sigma_\infty] = -\int_{\Sigma_0} [\Sigma_0] = 1$. We look for K\"ahler metrics whose K\"ahler classes have the form $b \lambda [\Sigma_\infty] - a \lambda [\Sigma_0]$ with $b > a > 0$. Note also that $c_1(L) = [-\sqrt{-1}\ddbar\rho] = -\lambda \pi^* [\omega_\Sigma] = \lambda [\Sigma_0] - \lambda [\Sigma_\infty]$.

In order to define a K\"ahler metric in the K\"ahler class $b [\Sigma_\infty] - a[\Sigma_0]$, we first define a momentum profile function $f(\rho)$ on $M_0 = M \backslash (\Sigma_0 \cup \Sigma_\infty)$. The idea of this momentum construction comes from the works \cite{HgSgr} by Hwang-Singer and \cite{ACGT} by Apostolov-Calderbank-Gauduchon-(T{\o}nnesen-Friedman) on extremal and constant scalar curvature K\"ahler metrics. Together with a pair of asymptotic conditions given below, one can extend the metric induced by $f$ to the whole manifold $M$. Here is the detail:

Let $f(\rho) : \R \to (a,b)$ be a strictly increasing function. We define a K\"ahler metric $\omega$ on $M_0$ by
\[\omega = f(\rho) \sqrt{-1}\ddbar{\rho} +\sqrt{-1} f_\rho(\rho) \partial \rho \wedge \partial \bar{\rho}.\]
\begin{remark}
If we let $u(\rho)$ be the anti-derivative of $f$, i.e. $u_\rho = f$, then one can check that $\omega = \sqrt{-1}\ddbar{u(\rho)}$ on $M_0$.
\end{remark}
In order for the K\"ahler metric to be defined on $M$, we require the following asymptotic conditions:
\begin{enumerate}
\item There exists a smooth function $F_0 : [0, \infty) \to \R$ with $F_0 (0) = a$ and $F_0'(0) > 0$, so that $f(\rho) = F_0 (e^{2\rho})$ as $\rho \to -\infty$.
\item There exists a smooth function $F_\infty : [0, \infty) \to \R$ with $F_\infty (0) = b$ and $F_\infty'(0) >0$ so that $f(\rho) = F_\infty (e^{-2\rho})$ as $\rho \to \infty$.
\end{enumerate}
Note that $f$ has to be a strictly increasing function, so we have $a < f(\rho) < b \text{ for } \rho \in \R$, and
\begin{align*}
\lim_{\rho \to -\infty} f(\rho) & = a,\\
\lim_{\rho \to \infty} f(\rho) & = b.\\
\lim_{\rho \to \pm \infty} f_\rho(\rho) & = 0.
\end{align*}
The K\"ahler class $[\omega]$ can be easily seen to be $[\omega] =b\lambda[\Sigma_\infty] - a\lambda[\Sigma_0]$, because
\begin{align*}
\langle [\omega], \Sigma_\infty \rangle & = \int_{\Sigma_\infty} b [\sqrt{-1}\ddbar{\rho}]\\
& = \int_{\Sigma_\infty} b\lambda[\Sigma_\infty]=  b\lambda,\\
\langle [\omega], \Sigma_0 \rangle & = \int_{\Sigma_0} a[\sqrt{-1}\ddbar{\rho}]\\
& = \int_{\Sigma_0} -a\lambda[\Sigma_0] = a\lambda.
\end{align*}

Under this construction, the K\"ahler form depends only on the height parameter $\rho$. We can see immediately that these K\"ahler metrics are invariant under the $U(1)$-action defined earlier, since the action preserves $\rho$: $\|e^{i\theta} u\|_h = \|u\|_h$ for any section $u \in \Gamma(\Sigma, L)$.

Note that for $(\Sigma, L) = (\mathbb{P}^n, \mathcal{O}_{\mathbb{P}^n}(-k))$, i.e. Hirzebruch-type manifolds, the above momentum construction with $\omega_\Sigma = \omega_{\text{FS}}$, i.e. the Fubini-Study metric, is the $U(n+1) / \Z_k$-symmetry initiated by Calabi in \cite{Cal82}.

Next we derive the local expression of the K\"ahler metric $\omega$ constructed by the above momentum profile as well as its Ricci curvature. Let $(z_1, \ldots, z_n, \xi)$ be local holomorphic coordinates of $M$ where $\mathbf{z} = (z_1, \ldots, z_n)$ are the base coordinates and $\xi$ is the fiber coordinate. Recall that the height parameter is defined to be
\[\rho = \log\|\cdot\|^2_h.\]
Let $\phi(\mathbf{z})$ be a positive function such that $\|\xi\|^2 = |\xi|^2 \phi(\mathbf{z})$ for any $(\mathbf{z}, \xi)$ in the local coordinate chart. Then we have
\begin{equation}\label{rho_coord}
\rho = \log|\xi|^2 + \log \phi(\mathbf{z}).
\end{equation}
Using this, one can easily check $\rho_{\xi\bar\xi} = \rho_{i\bar\xi} = \rho_{\xi\bar i} = 0$ for any $i = 1, 2, \ldots, n$.
Moreover, $\sqrt{-1}\ddbar\rho = \lambda \pi^*\omega_\Sigma$, so we can let $\lambda \pi^*\omega_\Sigma = \rho_{i\bar j} dz^i \wedge d\bar{z}^j$. Hence, the K\"ahler metric in $(\mathbf{z}, \xi)$ coordinates is given by
\begin{align*}
\omega & = \sqrt{-1} \sum_{i,j=1}^n (f \rho_{i\bar{j}} + f_\rho \rho_i \rho_{\bar{j}}) dz^i \wedge dz^{\bar{j}} + \sqrt{-1} f_\rho \sum_{i=1}^n \rho_i \rho_{\bar\xi} dz^i \wedge d\bar\xi\\
& \quad + \sqrt{-1} f_\rho \sum_{i=1}^n \rho_\xi \rho_{\bar i} d\xi \wedge dz^{\bar i} + \sqrt{-1} f_\rho |\rho_\xi|^2 d\xi \wedge d\bar\xi.
\end{align*}
Let $g$ be the metric associated to the K\"ahler form $\omega$, and $g_\Sigma$ be that of $\omega_\Sigma$. The determinant of the metric $g$ and its logarithm are given by
\begin{align*}
\det(g) & = \lambda^n f^n f_\rho \det(g_\Sigma) |\xi|^{-2},\label{det}\\
\log\det(g) & = n\log\lambda + n\log f + \log f_\rho + \log \det(g_\Sigma) - \log|\xi|^2.
\end{align*}
Using this, one can then compute the Ricci form $-\sqrt{-1}\ddbar\log\det(g)$:
\begin{align*}
\Ric(\omega) & = -\sqrt{-1}\ddbar\log\det(g)\\
& = \{(\nu \lambda^{-1} - \partial_\rho (n\log f + \log f_\rho))\rho_{i\bar j} - \partial_{\rho\rho} (n\log f + \log f_\rho) \rho_i \rho_{\bar j}\} dz^i \wedge dz^{\bar j}\\
& \quad - \partial_{\rho\rho} (n\log f + \log f_\rho) \rho_i \rho_{\bar \xi} dz^i \wedge d\bar\xi - \partial_{\rho\rho} (n \log f + \log f_\rho)\rho_\xi \rho_{\bar i} d\xi \wedge dz^{\bar i}\\
& \quad - \partial_{\rho\rho} (n\log f + \log f_\rho) |\rho_\xi|^2 d\xi \wedge d\bar\xi.
\end{align*}
In the computation of the Ricci form, we used the fact that $\omega_\Sigma$ is K\"ahler-Einstein so that $- \sqrt{-1}\ddbar\log\det(g_\Sigma) = \nu\omega_\Sigma$.

Observing that the $\omega$ and $\Ric$ have similar linear-algebraic expressions when $\omega$ is constructed by a momentum profile $f$, one can see easily that the K\"ahler-Ricci flow on $M$ is equivalent to a parabolic equation that evolves the momentum profile. In other words, the K\"ahler-Ricci flow preserves the momentum construction. Precisely, we have
\begin{proposition}
Suppose $\omega_0$ is the initial K\"ahler form on $M$ with momentum profile $f_0(\rho)$, then the solution $\omega_t,\, t \in [0,T)$ to the K\"ahler-Ricci flow $\partial_t \omega_t = -\Ric(\omega_t)$ also admits a momentum profile $f(\rho, t)$ at each time $t \in [0,T)$ where $f(\rho, t)$ evolves by
\begin{equation*}
\D{f}{t} = \dd{\rho}\left(n \log f(\rho, t) + \log f_\rho(\rho, t)\right) - \frac{\nu}{\lambda}, \quad f(\rho, 0) = f_0(\rho);
\end{equation*}
or equivalently,
\begin{equation}\label{eq_f}
\D{f}{t} = \frac{f_{\rho\rho}}{f_\rho} + n \frac{f_\rho}{f} - \frac{\nu}{\lambda}, \quad f(\rho, 0) = f_0(\rho).
\end{equation}
\end{proposition}

\section{K\"ahler classes under K\"ahler-Ricci Flow}
From now on, we will consider the K\"ahler-Ricci flow $\partial_t \omega_t = -\Ric(\omega_t)$ on $M$ which satisfies the aforesaid $U(1)$-symmetry and admits evolving momenta $f(\rho, t)$. We say $T$ is the blow-up time of the Ricci flow if $[0, T)$ is the maximal time interval for the Ricci flow to exist. For Ricci flow on compact K\"ahler manifolds, the blow-up time is completely determined by the initial K\"ahler class and the first Chern class. Namely, we have the following theorem proved by Tian-Zhang:
\begin{theorem}[Tian-Zhang, \cite{TnZg06}]\label{TZ}
Let $(X, \omega(t))$ be an (unnormalized) K\"ahler-Ricci flow $\partial_t \omega_t = -\Ric(\omega_t)$ on a compact K\"ahler manifold $X^n$. Then the blow-up time $T$ is given by
\begin{equation*}\label{maximal_time}
T = \sup \{ t : [\omega_0] + t c_1 (K_X) > 0\},
\end{equation*}
where $K_X := \det(T^*X)$ is the canonical line bundle of $X^n$.
\end{theorem}

Note that the K\"ahler class $[\omega_t]$ at any time $t$ is given by $[\omega_t] = [\omega_0] + t c_1(K_X)$. In order to work out the evolving K\"ahler classes and the blow-up time, one needs to understand the first Chern class of $K_X$, which can be computed by the adjunction formula.

Given any smooth divisor $D$ of compact K\"ahler manifold $X$, the adjunction formula relates $K_X$ and $K_D$ by
\begin{equation}\label{adj}
K_D = \eval{(K_X \otimes \mathcal{N}_{M/D})}{D},
\end{equation}
where $\mathcal{N}_{M/D}$ is the normal bundle of $D$ in $M$.

Using \eqref{adj}, one can easily work out $c_1(K_M)$ by taking $D = \Sigma_0, \Sigma_\infty$ in turn. For example, taking $D = \Sigma_\infty$, we have
\begin{align*}
K_{\Sigma_\infty} & = \eval{(K_M \otimes L^*)}{\Sigma_\infty},\\
\langle c_1(K_{\Sigma_\infty}), [\Sigma_\infty] \rangle & = \langle c_1(K_M) - c_1(L), [\Sigma_\infty] \rangle.
\end{align*}
Since $\Sigma$ is K\"ahler-Einstein such that $\Ric(\omega_\Sigma) = \nu\omega_\Sigma$, we then have
\begin{equation*}
\langle c_1(K_{\Sigma_\infty}), [\Sigma_\infty] \rangle = -\nu.
\end{equation*}
Since $c_1(L) = \lambda[\Sigma_0] - \lambda[\Sigma_\infty]$, we have $\langle c_1 (L), [\Sigma_\infty] \rangle = -\lambda$, and hence
\begin{equation*}
\langle c_1(K_M), [\Sigma_\infty] \rangle = -\nu - \lambda.
\end{equation*}
Similarly, one can also show by taking $D = \Sigma_0$ in \eqref{adj} (now $\mathcal{N}_{M \backslash D} = L$) to show
\begin{equation*}
\langle c_1(K_M), [\Sigma_0] \rangle = -\nu + \lambda.
\end{equation*}
Therefore, the first Chern class of the canonical line bundle $K_M$ is given by:
\begin{equation*}
c_1(K_M) = (-\nu - \lambda)[\Sigma_\infty] - (-\nu + \lambda)[\Sigma_0].
\end{equation*}
Hence, under the K\"ahler-Ricci flow $\partial_t \omega_t = -\Ric(\omega_t)$ with initial class $[\omega_0] = b_0 \lambda [\Sigma_\infty] - a_0 \lambda [\Sigma_0]$, the K\"ahler class evolves by
\begin{equation}\label{kahler_classes}
[\omega_t] = (b_0 \lambda - (\nu+\lambda)t) [\Sigma_\infty] - (a_0 \lambda - (\nu-\lambda)t) [\Sigma_0].
\end{equation}
We denote $[\omega_t] = \lambda b_t [\Sigma_\infty] - \lambda a_t [\Sigma_0]$ where $a_t, b_t$ are defined by
\begin{align}
a_t := a_0 - \frac{(\nu-\lambda)}{\lambda}t,\label{a_t}\\
b_t := b_0 - \frac{(\nu+\lambda)}{\lambda}t.\label{b_t}
\end{align}
Note also that $[\pi^*\omega_\Sigma] = [\Sigma_\infty] - [\Sigma_0]$, therefore the K\"ahler class can also be expressed as
\begin{equation}
[\omega_t] = \lambda a_t [\pi^*\omega_\Sigma] + \lambda(b_t - a_t)[\Sigma_\infty].
\end{equation}
Hence, by Theorem \ref{maximal_time}, the maximal time is characterized by $\lambda$ and $\nu$ in the following way:
\label{cases}
\begin{itemize}
\item Case 1: $\nu \leq \lambda$\\
In this case, $[\omega_t]$ ceases to be K\"ahler when $b_t = a_t$, namely, at $T := \frac{b_0 - a_0}{2}$. The limiting K\"ahler class is given by
\[[\omega_T] = \lambda a_T[\pi^*\omega_\Sigma].\]
This holds true for any given $b_0 > a_0 > 0$.
\item Case 2: $\nu > \lambda$\\
We further divide it into three sub-cases
\begin{enumerate}[(i)]
\item \underline{$(\nu - \lambda) b_0 < (\nu + \lambda) a_0$}:\\
$[\omega_t]$ ceases to be K\"ahler when $b_t = a_t$. Likewise, the limiting K\"ahler class is given by
\[[\omega_T] = \lambda a_T[\pi^*\omega_\Sigma].\]
\item \underline{$(\nu - \lambda) b_0 = (\nu + \lambda) a_0$}: \\
$[\omega_t]$ is then proportional to $c_1(K_M^{-1})$, i.e. canonical class. The flow stops at $T = \frac{a_0\lambda}{\nu-\lambda}$ and the limiting class $[\omega_T] = 0$. It is well-known (see e.g. \cite{SsTn08}) that in such case $(M, g(t))$ extincts and converges to a point in the Gromov-Hausdorff sense as $t \to T$.
\item \underline{$(\nu - \lambda) b_0 > (\nu + \lambda) a_0$}: \\
$[\omega_t]$ ceases to be K\"ahler when at $T = a_0$, and the limit class is given by $[\omega_T] = \lambda b_T[\Sigma_\infty]$.
\end{enumerate}
\end{itemize}

This trichotomy  resembles that in Song-Weinkove's work \cite{SgWk09} on Hirzebruch surfaces and Hirzebruch-type manifolds, i.e. $(\Sigma, L) = (\mathbb{P}^n, \mathcal{O}_{\mathbb{P}^n} \oplus \mathcal{O}_{\mathbb{P}^n}(-k))$. In their work, from which our study was motivated, similar trichotomy of the blow-up time of the K\"ahler-Ricci flow with initial K\"ahler class $[\omega_0]$ was also exhibited as it is characterized by the triple $(n, k, [\omega_0])$. It was shown in \cite{SgWk09} assuming Calabi's $U(n+1) / \Z_k$-symmetry and in \cite{SSW11} assuming $\Sigma$ is projective that in case of having limiting K\"ahler class $a_T[\pi^*\omega_\Sigma]$, the K\"ahler-Ricci flow collapses the $\mathbb{P}^1$-fiber of the projective bundle, which hereof converges to some K\"ahler metric of $\Sigma$ as metric spaces in Gromov-Hausdorff sense.

Case 2(iii) is in reminiscence of Song-Weinkove's recent works \cite{SgWk10} and \cite{SgWk11} of contracting exceptional divisors, in which $\mathcal{O}(-k)$-blow-up of arbitrary compact K\"ahler manifold $X$ are considered. In their works, a cohomological condition is given on the initial K\"ahler class and the first Chern class, under which the blown-up manifold will converge in Gromov-Hausdorff sense back to $X$ with orbifold singularity of type $\mathcal{O}(-k)$. There is no symmetry assumption in these works.

In our paper, we will only focus on Case 1 and Case 2(i) which exhibit collapsing of $\mathbb{P}^1$-fiber assuming the K\"ahler metric admits the aforesaid momentum construction.

\section{Estimates of the K\"ahler-Ricci flow}
From now on we assume that the triple $(\Sigma, L, [\omega_0])$ satisfies Case 1 or Case 2(i) stated in the previous section, i.e. either
\begin{align*}
\nu & \leq \lambda, \text{ or }\\
\nu & > \lambda \text{ and } (\nu - \lambda)b_0 < (\nu + \lambda)a_0. 
\end{align*}
Recall that $\nu$ is the Ricci curvature of the K\"ahler-Einstein manifold $\Sigma$ and $\lambda$ is the Chern curvature of the Hermitian-Einstein line bundle $L$, i.e.
\begin{align*}
\Ric(\omega_\Sigma) & = \nu\omega_\Sigma,\\
\sqrt{-1}\ddbar\rho & = -\lambda \pi^*\omega_\Sigma.
\end{align*}
Recall that the first Chern class of $K_M$ and the evolving K\"ahler class are given by:
\begin{align*}
c_1(K_M) & = (-\nu - \lambda)[\Sigma_\infty] - (-\nu + \lambda)[\Sigma_0]\\
& = (-\nu+\lambda)[\pi^*\omega_\Sigma] - 2\lambda[\Sigma_\infty],\\
[\omega_t] & = \lambda b_t [\Sigma_\infty] - \lambda a_t [\Sigma_0]\\
& = \lambda a_t [\pi^*\omega_\Sigma] + 2\lambda(T-t) [\Sigma_\infty]
\end{align*}
where $a_t$ and $b_t$ defined in \eqref{a_t} and \eqref{b_t}.

Since pluripotential theory plays a very important role in K\"ahler-Ricci flow and in K\"ahler geometry in general, we would like understand the K\"ahler-Ricci flow $\partial_t \omega_t = -\Ric(\omega_t)$ from potential functions viewpoint. To do so, we need a reference family of K\"ahler metrics $\{\hat{\omega}_t\}_{t \in [0, T)}$ whose K\"ahler class at each time $t$ coincides with that of $\omega_t$, the K\"ahler-Ricci flow solution. We choose $\hat{\omega}_t$ to be the $U(1)$-invariant K\"ahler metric induced by the following momentum profile:
\begin{equation*}
\hat{f}(\rho, t) := a_t +\frac{(b_t - a_t)e^{2\rho}}{1+e^{2\rho}} = a_t + \frac{2\lambda(T-t)e^{2\rho}}{1+e^{2\rho}}.
\end{equation*}
This momentum profile gives the following K\"ahler metric:
\begin{equation*}
\hat{\omega}_t = a_t \sqrt{-1}\ddbar\rho + 2\sqrt{-1}\lambda(T-t) \left(\frac{e^{2\rho}}{1+e^{2\rho}} \ddbar\rho + \frac{2e^{2\rho}}{(1+e^{2\rho})^2} \partial \rho \wedge \bar{\partial} \rho \right).
\end{equation*}
Clearly, $\hat{f}$ satisfies the asymptotic conditions for extending $\hat{\omega}_t$ to the whole $M$. Also, we have $[\hat{\omega}_t] = [\omega_t]$ because $\hat{f} \to a_t$ as $\rho \to -\infty$ and $\hat{f} \to b_t$ as $\rho \to \infty$.

For simplicity, we denote $\Theta := \frac{e^{2\rho}}{1+e^{2\rho}} \sqrt{-1}\ddbar\rho + \frac{2e^{2\rho}}{(1+e^{2\rho})^2} \sqrt{-1} \partial \rho \wedge \bar{\partial} \rho$, so that
\begin{equation*}
\hat{\omega}_t = a_t \pi^*\omega_\Sigma + 2\lambda(T-t) \Theta.
\end{equation*}
Note that $[\Theta] = [\Sigma_\infty]$ and so $\D{\hat{\omega}_t}{t}=(-\nu + \lambda)\pi^*\omega_\Sigma - 2\lambda\Theta \in c_1(K_M)$. Take $\Omega$ be a fixed volume form of $M$ such that $\D{\hat{\omega}_t}{t} = \sqrt{-1}\ddbar\log\Omega$. Then the K\"ahler-Ricci flow $\partial_t \omega_t = -\Ric(\omega_t)$ is equivalent to the following complex Monge-Amp\`ere equation
\begin{equation}\label{mgap}
\D{\phi}{t} = \log \frac{\det(\hat{\omega}_t + \sqrt{-1}\ddbar\phi)}{(T-t)\Omega}, \quad \eval{\phi}{t=0} = \phi_0
\end{equation}
in a sense that $\omega_t = \hat{\omega}_t + \sqrt{-1}\ddbar\phi$, $t \in [0,T)$ is a solution to the K\"ahler-Ricci flow $\partial_t \omega_t = -\Ric(\omega_t)$ with initial data $\omega_0 = \hat{\omega}_0 + \sqrt{-1}\ddbar\phi_0$ if and only if $\phi : M \times [0,T)$ is a solution to \eqref{mgap}.

Working similarly as in \cite{SgTn06, SgWk09, SgWk10, TnZg06} etc., one can derive the following estimates using maximum principles.
\begin{lemma}\label{est_1}
There exists a constant $C = C(n, \omega_0, \nu, \lambda)> 1$ such that the following holds
\begin{enumerate}
\item $|\phi(t)| \leq C$,
\item $\omega_t^{n+1} \leq C \Omega$, and
\item $\Tr_{\omega_t} \pi^*\omega_\Sigma \leq C$.
\end{enumerate}
\end{lemma}
\begin{proof}
The proof goes similarly as in \cite{SgWk09}. First note that since
\[\hat{\omega}_t^{n+1} \geq 2\lambda a_t (n+1)(T-t) (\pi^*\omega_\Sigma)^n \wedge \Theta,\]
one can then find constant $C > 0$ independent of $t$ such that 
\begin{equation}\label{omega_hat_bd}
C^{-1} (T-t)\Omega \leq \hat{\omega}^{n+1}_t \leq C(T-t)\Omega.
\end{equation}
Consider the function $\tilde{\phi} = \phi + (1+\log C)t$, at the point $p_t \in M$ where $\tilde{\phi}$ achieves its minimum at time $t$, we have $\ddbar\tilde{\phi} = \ddbar\phi \geq 0$. Therefore,
\begin{align*}
\ddt \tilde{\phi}_{\min}(t) & = \eval{\log\frac{\det(\hat{\omega}_t + \sqrt{-1}\ddbar\phi)}{(T-t)\Omega}}{p_t} + \log C\\
& \geq \eval{\log\frac{\det{\hat{\omega}_t}}{C^{-1}(T-t)\Omega}}{p_t} \geq 0.
\end{align*}
Here we used \eqref{omega_hat_bd}. It proves $\phi$ is uniformly bounded from below as the flow encounters finite-time singularity. The uniform upper bound for $\phi$ follows similarly. 

For (2), we consider $Q := \D{\phi}{t} - |\lambda - \nu|a^{-1}\phi + \log(T-t)$ where $a := \inf_{[0,T)} a_t > 0$. By direct computation, we have
\begin{align}
\D{Q}{t} & = \Tr_{\omega_t}((\lambda-\nu) \pi^*\omega_\Sigma - 2\lambda\Theta) + \Delta \left(\D{\phi}{t}\right)  \label{ddt_Q}\\
& \quad -|\lambda-\nu|a^{-1}(Q+|\lambda-\nu|a^{-1}\phi - \log(T-t)) - \frac{1}{T-t} \notag\\
& = \Delta Q + |\lambda-\nu|a^{-1} \Delta\phi + (\lambda-\nu) \Tr_{\omega_t}\pi^*\omega_\Sigma - 2\lambda\Tr_{\omega_t}\Theta \notag\\
& \quad -|\lambda-\nu|a^{-1}(Q+|\lambda-\nu|a^{-1}\phi - \log(T-t)) - \frac{1}{T-t}. \notag
\end{align}
Since $\omega_t = a_t \pi^*\omega_\Sigma + 2\lambda(T-t)\Theta + \sqrt{-1}\ddbar\phi(t)$, taking trace with respect to $\omega_t$ yields
\[n+1 = a_t \Tr_{\omega_t}\pi^*\omega_\Sigma + 2\lambda(T-t)\Tr_{\omega_t}\Theta + \Delta\phi \geq a_t \Tr_{\omega_t}\pi^*\omega_\Sigma + \Delta\phi.\]
Hence we have, 
\begin{equation}\label{lap_phi}
|\lambda-\nu|a^{-1}\Delta\phi \leq |\lambda-\nu|a^{-1}(n+1) - |\lambda-\nu|\Tr_{\omega_t}\pi^*\omega_\Sigma.
\end{equation}
Note that $a_t \geq a$ for any $t \in [0,T)$. Combining \eqref{ddt_Q} and \eqref{lap_phi}, we have
\begin{equation}\label{box_Q}
\square Q \leq (n+1)|\lambda-\nu|a^{-1} - |\lambda-\nu|a^{-1} (Q + |\lambda-\nu|a^{-1}\phi) + |\lambda-\nu|a^{-1}\log T.
\end{equation}
As $\phi$ is uniformly bounded from (1), \eqref{box_Q} implies a uniformly upper bound for $Q$. Since $Q = \log\frac{\det \omega_t}{\Omega} - |\lambda-\nu|a^{-1}\phi$, again together with the uniform bound for $\phi$, we proved (2).

Finally, for (3), we let $(z_1, \ldots, z_n, \xi)$ be local holomorphic coordinates such that $\mathbf{z}=(z_1,\ldots,z_n)$ is the base coordinate and $\xi$ is the fiber coordinate. Then the bundle map is given by $\pi : (\mathbf{z}, \xi) \mapsto \mathbf{z}$. Write $\lambda \pi^*\omega_\Sigma = \rho_{i\bar j} dz^i \wedge d\bar{z}^j$. Assume the K\"ahler metrics $\omega_t$ admits momentum profiles $f(\rho, t)$, we have $g^{i\bar j} = \frac{1}{f}\rho^{i\bar j}$ and one can prove
\begin{equation*}
\Tr_{\omega_t}\pi^*\omega_\Sigma  = \frac{1}{\lambda} g^{i\bar j} \rho_{i\bar j} = \frac{n}{\lambda f}
\end{equation*} which is clearly bounded from above uniformly independent of $t$.
\end{proof}
Next, we will derive estimates on the K\"ahler-Ricci flow by assuming the metric is $U(1)$-invariant and admits a momentum profile $f(\rho, t)$. First note that because $f_\rho(\rho, t) > 0$ for any $t$ and also $\lim_{\rho\to -\infty}f(\rho, t) = a_t$, $\lim_{\rho\to\infty}f(\rho,t) = b_t$, we have
\begin{equation*}
a_t < f(\rho, t) < b_t, \quad \text{ for any } (\rho, t) \in \R \times [0,T).
\end{equation*}
Note that $a_t$ and $b_t$ are both bounded away from zero as $t \to T$, (2) in Lemma \ref{est_1} implies $f_\rho$ is also uniformly bounded. Using these, one is able to derive the following estimates.
\begin{lemma}\label{est_2}
There exists a constant $C = C(n, \omega_0, \nu, \lambda) > 0$ such that
\begin{enumerate}
\item $C^{-1} \leq f \leq C$,
\item $\abs{\frac{f_{\rho\rho}}{f_\rho}} \leq C$,
\item $f_{\rho} \leq C(T-t)$
\end{enumerate}
for any $(\rho, t) \in \R \times [0,T)$.
\end{lemma}
\begin{proof}
As discussed above, (1) clearly holds because $a_t$ is bounded away from zero and $b_t$ is uniformly bounded above on $[0, T)$.

For (2), first note that by the asymptotic conditions of the momentum profile $f(\rho, t)$, $\lim_{\rho\to\pm\infty} \abs{\frac{f_{\rho\rho}}{f_\rho}} = 2$ for any $t \in [0,T)$, so $\sup_{\R \times [0,T-\epsilon)} \abs{\frac{f_{\rho\rho}}{f_\rho}}$ exists for every $\epsilon > 0$. We will derive the uniform lower bound for $\frac{f_{\rho\rho}}{f_\rho}$ on $[0,T)$ since the upper bound is similar. Given any $\epsilon > 0$, let $(\rho_{\epsilon}, t_{\epsilon}) \in \R \times [0,T-\epsilon)$ be the point such that
\[\eval{\frac{f_{\rho\rho}}{f_\rho}}{(\rho_{\epsilon}, t_{\epsilon})} = \sup_{\R \times [0,T-\epsilon)} \frac{f_{\rho\rho}}{f_\rho}.\]
Then at $(\rho_{\epsilon}, t_{\epsilon})$, one has $\dd{t}\left(\frac{f_{\rho\rho}}{f_\rho}\right) \geq 0$, $\dd{\rho}\left(\frac{f_{\rho\rho}}{f_\rho}\right) = 0$, and $\frac{\partial^2}{\partial \rho^2}\left(\frac{f_{\rho\rho}}{f_\rho}\right) \leq 0$.

Recall that $f$ satisfies heat-type equation \eqref{eq_f}, i.e. $\D{f}{t} = \frac{f_{\rho\rho}}{f_\rho} + n\frac{f_\rho}{f} - \frac{\nu}{\lambda}$. By direct computation, one has
\begin{align*}
\dd{t}\left(\frac{f_{\rho\rho}}{f_\rho}\right) & = \frac{2 n f_\rho^2}{f^3}-\frac{2 n f_{\rho\rho}}{f^2}-\frac{n f_{\rho\rho}^2}{f f_\rho^2}+\frac{3 f_{\rho\rho}^3}{f_\rho^4}+\frac{n f_{\rho\rho\rho}}{f f_\rho}-\frac{4 f_{\rho\rho} f_{\rho\rho\rho}}{f_\rho^3}+\frac{f_{\rho\rho\rho\rho}}{f_\rho^2},\\
\dd{\rho}\left(\frac{f_{\rho\rho}}{f_\rho}\right) &= \frac{f_\rho f_{\rho\rho\rho} - f_{\rho\rho}^2}{f_\rho^2},\\
\frac{\partial^2}{\partial \rho^2}\left(\frac{f_{\rho\rho}}{f_\rho}\right) & = \frac{2f_{\rho\rho}^3}{f_\rho^3} - \frac{3f_{\rho\rho}f_{\rho\rho\rho}}{f_\rho^2} + \frac{f_{\rho\rho\rho\rho}}{f_\rho}.
\end{align*}
Evaluating at $(\rho_\epsilon, t_\epsilon)$, the fact that $\dd{\rho}\left(\frac{f_{\rho\rho}}{f_\rho}\right) = 0$ implies $f_{\rho\rho\rho} = \frac{f^2_{\rho\rho}}{f_\rho}$ at $(\rho_\epsilon, t_\epsilon)$. By substituting $f_{\rho\rho\rho} = \frac{f^2_{\rho\rho}}{f_\rho}$ into the expressions of $\dd{t}\left(\frac{f_{\rho\rho}}{f_\rho}\right)$ and $\frac{\partial^2}{\partial \rho^2}\left(\frac{f_{\rho\rho}}{f_\rho}\right)$, one can check that after cancellation of terms, we have
\begin{equation*}
0 \leq \left(\dd{t} - \frac{1}{f_\rho} \frac{\partial^2}{\partial \rho^2} \right) \frac{f_{\rho\rho}}{f_\rho} = \frac{2nf_\rho^2}{f^3} - \frac{2nf_{\rho\rho}}{f^2} \quad \text{ at } (\rho_\epsilon, t_\epsilon).
\end{equation*}
It shows $\sup_{\R \times [0, T-\epsilon)} \frac{f_{\rho\rho}}{f_\rho} = \eval{\frac{f_{\rho\rho}}{f_\rho}}{(\rho_\epsilon, t_\epsilon)} \leq \eval{\frac{f_\rho}{f}}{(\rho_\epsilon, t_\epsilon)}$. Since $f_\rho$ is uniformly bounded from above and $f > C^{-1}$, there exists $C > 0$ independent of $\epsilon$ such that $\sup_{\R \times [0, T-\epsilon)} \frac{f_{\rho\rho}}{f_\rho} \leq C$. Similar approach proves $\inf_{\R \times [0,T)} \frac{f_{\rho\rho}}{f_\rho} \geq -\tilde{C}$ for some uniform constant $\tilde{C} > 0$. It completes the proof of (2).

Part (3) follows from part (2). Precisely, (2) implies $|(\log f_\rho)_\rho| \leq C$. If we let $\rho_t \in \R$ such that $f_\rho(\rho_t) = \sup_{\rho\in\R} f_\rho$. Then by the mean-value theorem,
\[|\log f_\rho(\rho, t) - \log f_\rho(\rho_t, t)| \leq C|\rho - \rho_t|.\]
Thus for $\rho \in [\rho_t - C^{-1}, \rho_t + C^{-1}]$, we have
\[\log \left(\frac{f_\rho(\rho,t)}{f_\rho(\rho_t,t)}\right) \geq -1,\]
or equivalently, $f_\rho(\rho, t) \geq e^{-1} f_\rho(\rho_t, t)$. We then have
\begin{equation*}
\int_{\R} f_\rho d\rho \geq \int_{\rho_t - C^{-1}}^{\rho_t + C^{-1}} f_\rho d\rho \geq 2C^{-1} e^{-1}f_\rho(\rho_t,t).
\end{equation*}
On the other hand, we have
\[\int_\R f_\rho d\rho = f(\infty) - f(-\infty) = b_t - a_t = 2\lambda(T-t).\]
Hence $\sup_{\rho\in\R} f_\rho \leq C(T-t)$ for some uniform constant $C$.
\end{proof}
Lemma \ref{est_2} implies the $\mathbb{P}^1$-fiber of our manifold $M$ is collapsing along the flow. Precisely we have the following:
\begin{proposition}\label{fiber_collapse}
Assume $(\Sigma, L, [\omega_0])$ satisfies the condition stated in Case 1 and Case 2(i) in P.\pageref{cases}. Let $V_x \in T_x M$ be a tangent vector of $M$ at $x \in M \backslash (\Sigma_0 \cup \Sigma_\infty)$ which lies $T_x \mathbb{P}_x^1$. Here we denote $\mathbb{P}_x^1$ as the $\mathbb{P}^1$-fiber passing through $x$. Then we have $\|V_x\|_{g(t)} \to 0$ as $t \to T$.
\end{proposition}
\begin{proof}
It suffices to express $\|V_x\|_{g(t)}$ in terms of $f$ and $f_\rho$. Since the metric $g(t)$ is given by
\[g(t) = f \lambda \pi^*g_\Sigma + f_\rho \partial\rho \otimes \bar\partial\rho.\]
Since $V_x$ is parallel to the fiber, we have $\pi_* V_x = 0$ and so $\pi^*g_\Sigma(V_x, \bar V_x) = 0$. Hence $\|V_x\|^2_{g(t)} = f_\rho \D{V_x}{\rho} \D{\bar V_x}{\rho} \to 0$ as $t \to T$. Here we have used part (3) of Lemma \ref{est_2}.
\end{proof}
Furthermore, Lemmas \ref{est_1} and \ref{est_2} provide enough estimates in order to show $(M, \omega_t)$ converges to $(\Sigma, a_T \omega_\Sigma)$ as metric spaces in Gromov-Hausdorff sense. 
\begin{theorem}\label{GrHdf}
Suppose $(\Sigma, L, [\omega_0])$ satisfies the condition stated in Case 1 and Case 2(i) in P.\pageref{cases}, then $(M, g(t))$ converges to the K\"ahler-Einstein manifold $(\Sigma, a_T \omega_\Sigma)$ in Gromov-Hausdorff sense as $t \to T$.
\end{theorem}
\begin{proof}
The proof goes almost the same as in Song-Weinkove's paper \cite{SgWk09} on Hirzebruch surfaces with Calabi ansatz. We will sketch the main idea here. For detail, please refer to Song-Weinkove's paper. The main ingredients of the argument are as follows:
\begin{enumerate}
\item the metric $g(t)$ is degenerating along the fiber direction on compact subsets of $M \backslash (\Sigma_0 \cup \Sigma_\infty)$,
\item $g(t)$ is bounded above uniformly $g(0)$, and
\item for any $0 < \alpha < 1$, $g(t)$ converges to $a_T \pi^*\omega_\Sigma$ in $C^{\alpha}$-sense on compact subsets of $M \backslash (\Sigma_0 \cup \Sigma_\infty)$.
\end{enumerate}
We have proved (1) in Proposition \ref{fiber_collapse}. (2) can be proved by a uniform estimate on $f_\rho$ which can be obtained easily by the bound on the volume form $\omega_t^{n+1}$ in Lemma \ref{est_1}. For (3), note that $\omega_t = f(\rho, t) \sqrt{-1}\ddbar\rho + \sqrt{-1} f_\rho(\rho, t) \partial\rho \wedge \bar\partial\rho$. One can compute that $\|\nabla_{g_0} g(t)\|^2_{g_0}$ is a polynomial expression of $f(\rho,t), f_\rho(\rho,t)$ and $f_{\rho\rho}(\rho,t)$ where the coefficients are time-independent. Lemma \ref{est_2} then shows for any compact subset $K \in M \backslash (\Sigma_0 \cup \Sigma_\infty)$, so we have $\sup_{K \times [0,T)} \|\nabla_{g_0} g(t)\|^2_{g_0} \leq C_K$ for some time independent constant $C_K > 0$. It proves (3).

To show the Gromov-Hausdorff convergence, first fix a leave of $\Sigma$ in $M \backslash (\Sigma_0 \cup \Sigma_\infty)$. We denote it by $\sigma(\Sigma)$. Using (2), one can choose a sufficiently small tubular neighborhood of $\Sigma_0$ and $\Sigma_\infty$ such that their complement contains $\sigma(\Sigma)$. Then given any two points $x_1, x_2 \in M$, we project them down to the base $\Sigma$ via the bundle map $\pi$. Consider the length of the geodesic $\gamma$ joining $\pi(x_1)$ and $\pi(x_2)$, by lifting the geodesic up by $\sigma$, we know that the lifted $\gamma$ has length arbitrarily close to the $a_T \omega_\Sigma$-length by (3). Finally, using (1), one can show $x_i$ is arbitrarily close to $\sigma \circ \pi (x_i)$ as $t \to T$. Using triangle inequality, one can then prove the $g(t)$-distance between $x_1$ and $x_2$ are is arbitrarily close to the $a_T\omega_\Sigma$-distance as $t \to T$.
\end{proof}
\section{Splitting Lemma}
In the singularity analysis of closed (real) 3-manifolds as in \cite{Ham95} and \cite{Perl1}, one often consider a rescaled dilation, which is a rescaled sequence of metrics $g_i (t) = K_i g(t_i + K^{-1}_i t)$ where $K_i$ are chosen such that $K_i = \|\Rm(x_i)\|_{g(t_i)} \to \infty$ and $\|\Rm_{g_i(t)}\|_{g_i(t)} \leq C$ for some uniform constant $C > 0$ independent of $i$. By Hamilton's compactness \cite{Ham95} and Perelman's local non-collapsing theorem \cite{Perl1}, one can extract a subsequence, still call it $g_i(t)$, such that $(M, g_i (t), x_i) \to (M_{\infty}, g_{\infty}(t), x_{\infty})$ on compact subsets in Cheeger-Gromov sense. The convergence is in $C^\infty$-topology because once the curvature tensor is uniformly bounded, Shi's derivative estimate in \cite{Shi} asserts all the higher order derivatives of $\Rm$ are uniformly bounded. The limit obtained is often called a singularity model. According to the curvature blow-up rate (Type I or II), a singularity model may be an ancient or eternal solution, and is $\kappa$-non-collapsed by Perelman's result. These singularity models encode crucial geometric data near the singularity region of the flow.

We will show that under our momentum construction and our assumption on the triple $(\Sigma, L, [\omega_0])$, the singularity model $M_{\infty}$ obtained by the aforesaid rescaled dilations splits isometrically into a product $N \times L$, where $\dim_\C N = n$ and $\dim_\C L = 1$. 

Let $(z_1, \ldots, z_n, \xi)$ be local holomorphic coordinates where $z = (z_1, \ldots, z_n)$ are the base coordinates and $\xi$ is the fiber coordinate. Then $\lambda \pi^*\omega_\Sigma = \sqrt{-1} \rho_{i\bar{j}}(z) dz^i \wedge dz^{\bar{j}}$, the the K\"ahler metric defined by momentum profile $f(\rho, t)$, its inverse and the Ricci tensor are locally written as
\begin{equation*}
g_{AB} = \begin{cases}
f\rho_{i\bar j} + f_\rho \rho_i \rho_{\bar j} & \text{ if } (A, B) = (i, \bar j)\\
f_\rho \rho_i \rho_{\bar\xi} & \text{ if } (A, B)= (i, \bar\xi)\\
f_\rho |\rho_\xi|^2 & \text{ if } (A, B) = (\xi, \bar\xi)
\end{cases},
\end{equation*}
\begin{equation*}
g^{AB} = \begin{cases}
\frac{1}{f} \rho^{i\bar j} & \text{ if } (A, B) = (i, \bar{j}) \\
- \frac{1}{f \rho_{\bar\xi}} \sum_{k=1}^n \rho^{i\bar k}\rho_{\bar k} & \text{ if } (A, B) = (i, \bar\xi)\\
\frac{1}{|\rho_\xi|^2}\left(\frac{1}{f_\rho} + \frac{\sum_{k,l=1}^n \rho^{k\bar l}\rho_k \rho_{\bar l}}{f}\right) & \text{ if } (A, B) = (\xi, \bar\xi)
\end{cases},
\end{equation*}
\begin{equation*}
\Ric_{AB} = \begin{cases}
(\nu\lambda^{-1} - F_\rho)\rho_{i\bar j} - F_{\rho\rho}\rho_i \rho_{\bar j} & \text{ if } (A, B) = (i, \bar j)\\
-F_{\rho\rho}\rho_i \rho_{\bar\xi} & \text{ if } (A, B) = (i, \bar\xi)\\
-F_{\rho\rho} |\rho_\xi|^2 & \text{ if } (A, B) = (\xi, \bar\xi)
\end{cases}
\end{equation*}
where $F = n\log f + \log f_\rho$.

From the local expressions of $g$ and $g^{-1}$, one can easily derive local expressions of the Christoffel symbols which we will need for deriving our splitting result.
\begin{lemma}\label{christoffel}
The Christoffel symbols of the K\"ahler metric $g$ on $M$ constructed by momentum profile $f$ are given by
\begin{align*}
\Gamma_{\xi\xi}^i & = 0,\\
\Gamma_{\xi\xi}^{\xi} & = \frac{f_{\rho\rho}}{f_\rho} \rho_\xi + \frac{\rho_{\xi\xi}}{\rho_\xi} = \left(\frac{f_{\rho\rho}}{f_\rho} - 1\right) \rho_\xi,\\
\Gamma_{i\xi}^\xi & = \left(\frac{f_{\rho\rho}}{f_\rho} - \frac{f_\rho}{f} \right) \rho_i,\\
\Gamma_{i\xi}^j & = \frac{f_\rho}{f} \delta_i^j \rho_{\xi},\\
\Gamma_{ij}^\xi & = \left(\frac{f_{\rho\rho}}{f_\rho} - \frac{2f_\rho}{f}\right) \frac{\rho_i\rho_j}{\rho_\xi} - \frac{1}{\rho_\xi}\left(\rho^{l\bar k}\rho_l\rho_{j\bar k i} + \rho_{ij}\right),\\
\Gamma_{ij}^k & = \frac{f_\rho}{f}(\rho_i \delta^k_j + \rho_j \delta^k_i) + \rho^{k\bar l}\rho_{j\bar l i}.
\end{align*}
\end{lemma}
\begin{remark}
Recall that for K\"ahler manifolds, the only (possibly) non-zero Christoffel symbols are those with indexes of either all $(1,0)$-type or all $(0,1)$-type. For succinctness, please excuse us for omitting those which are vanishing or conjugate to one of the above.
\end{remark}
\begin{remark}
We will see that the vanishing of $\Gamma_{\xi\xi}^i$ is crucial when dealing with the curvature tensor in the blow-up analysis in the next section. Moreover, we only need the first four Christoffel symbols in order to obtain the splitting lemma.
\end{remark}
\begin{proof}
Using the formula $\Gamma_{\alpha\beta}^\gamma = g^{\gamma\bar\delta}\partial_\alpha g_{\beta\bar\delta}$ for K\"ahler manifolds, one can compute the Christoffel symbols directly:
\begin{align*}
\Gamma_{\xi\xi}^i & = g^{i\bar j} \dd{\xi} g_{\xi\bar j} + g^{i\bar\xi} \dd{\xi} g_{\xi\bar\xi}\\
& = \frac{1}{f} \rho^{i\bar j} \dd{\xi}(f_\rho\rho_\xi\rho_{\bar j}) - \frac{1}{f\rho_{\bar\xi}}\rho^{i\bar k}\rho_{\bar k} \dd{\xi}(f_\rho\rho_\xi\rho_{\bar\xi})\\
& = \frac{1}{f}\rho^{i\bar j} (f_{\rho\rho}\rho_\xi\rho_\xi\rho_{\bar j} + f_\rho \rho_{\xi\xi} \rho_{\bar j}) - \frac{1}{f\rho_{\bar\xi}} \rho^{i\bar k}\rho_{\bar k} (f_{\rho\rho} \rho_\xi \rho_\xi \rho_{\bar\xi} + f_\rho\rho_{\xi\xi}\rho_{\bar\xi}) \\
& = 0.\\
\Gamma_{\xi\xi}^{\xi} & = \sum_{i=1}^n g^{\xi\bar i} \dd{\xi} g_{\xi\bar i} + g^{\xi\bar\xi}\dd{\xi} g_{\xi\bar\xi}\\
& = -\frac{1}{f\rho_\xi} \sum_{k=1}^n \rho^{k\bar i} \rho_k \dd{\xi} (f_\rho \rho_\xi \rho_{\bar i}) + \frac{1}{|\rho_\xi|^2} \left(\frac{1}{f_\rho} + \frac{\sum_{k,l=1}^n \rho^{k\bar l}\rho_k \rho_{\bar l}}{f}\right) \dd{\xi} (f_\rho \rho_\xi \rho_{\bar\xi})\\
& = -\frac{1}{f\rho_\xi} \sum_{k=1}^n \rho^{k\bar i} \rho_k \dd{\xi} (f_{\rho\rho} \rho_\xi^2 \rho_{\bar i} + f_\rho \rho_{\xi\xi} \rho_{\bar i})\\
& \quad + \frac{1}{|\rho_\xi|^2}\left(\frac{1}{f_\rho} + \frac{\sum_{k,l=1}^n \rho^{k\bar l}\rho_k \rho_{\bar l}}{f}\right) (f_{\rho\rho} \rho_\xi^2 \rho_{\bar\xi} + f_\rho \rho_{\xi\xi} \rho_{\bar\xi})\\
& = \frac{f_{\rho\rho}}{f_\rho} \rho_\xi + \frac{\rho_{\xi\xi}}{\rho_\xi} = \left(\frac{f_{\rho\rho}}{f_\rho} - 1\right) \rho_\xi.\\
\Gamma_{i\xi}^\xi & = \sum_{j=1}^n g^{\xi\bar j} \dd{z_i} g_{\xi\bar j} + g^{\xi\bar\xi} \dd{z_i} g_{\xi\bar\xi}\\
& = -\frac{1}{f\rho_\xi} \rho^{k\bar j} \rho_k \dd{z_i} (f_\rho \rho_\xi \rho_{\bar j}) + \frac{1}{|\rho_\xi|^2}\left(\frac{1}{f_\rho} + \frac{\sum_{k,l=1}^n \rho^{k\bar l}\rho_k \rho_{\bar l}}{f}\right) \dd{z_i} (f_\rho \rho_\xi \rho_{\bar \xi})\\
& = -\frac{1}{f\rho_\xi} \rho^{k\bar j} \rho_k (f_{\rho\rho} \rho_i \rho_\xi \rho_{\bar j} + f_\rho \rho_\xi \rho_{i\bar j}) \\
& \quad + \frac{1}{|\rho_\xi|^2}\left(\frac{1}{f_\rho} + \frac{\sum_{k,l=1}^n \rho^{k\bar l}\rho_k \rho_{\bar l}}{f}\right) (f_{\rho\rho} \rho_i \rho_\xi \rho_{\bar\xi})\\
& = \left(\frac{f_{\rho\rho}}{f_\rho} - \frac{f_\rho}{f} \right) \rho_i.\\
\Gamma_{i\xi}^j & = \sum_{k=1}^n g^{j\bar k} \dd{z_i} g_{\xi\bar k} + g^{j\xi} \dd{z_i} g_{\xi\bar\xi}\\
& = \frac{1}{f} \rho^{j\bar k} \dd{z_i} (f_\rho \rho_\xi \rho_{\bar k}) - \frac{1}{f\rho_{\bar\xi}} \rho^{j\bar k} \rho_{\bar k} \dd{z_i} (f_\rho \rho_\xi \rho_{\bar \xi})\\
& = \frac{1}{f} \rho^{j\bar k} (f_{\rho\rho} \rho_i \rho_\xi \rho_{\bar k} + f_\rho \rho_\xi \rho_{i\bar k}) - \frac{1}{f\rho_{\bar\xi}} \rho^{j\bar k} \rho_{\bar k} (f_{\rho\rho} \rho_i \rho_\xi \rho_{\bar\xi})\\
& = \frac{f_\rho}{f} \delta_i^j \rho_{\xi}.
\end{align*}
\begin{align*}
\Gamma_{ij}^\xi & = \sum_{k=1}^n g^{\xi\bar k}\dd{z_i} g_{j\bar k} + g^{\xi\bar\xi} \dd{z_i} g_{j\bar\xi}\\
& = -\frac{1}{f\rho_\xi}\rho^{\bar k l}\rho_l \dd{z_i}(f\rho_{j\bar k} + f_\rho \rho_j \rho_{\bar k}) + \frac{1}{|\rho_\xi|^2} \left(\frac{1}{f_\rho} + \frac{\rho^{k\bar l}\rho_k\rho_{\bar l}}{f}\right)\dd{z_i}(f_\rho\rho_j\rho_{\bar\xi})\\
& = -\frac{f_\rho}{f\rho_\xi}\rho^{l\bar k}\rho_l (\rho_i \rho_{j\bar k} + \rho_j\rho_{i\bar k}) - \frac{1}{\rho_\xi} \rho^{l\bar k}\rho_{j\bar ki}\rho_l + \frac{f_{\rho\rho}}{f_\rho \rho_\xi} \rho_i \rho_{\bar j} + \frac{1}{\rho_\xi}\rho_{i\bar j}\\
& = \left(\frac{f_{\rho\rho}}{f_\rho} - \frac{2f_\rho}{f}\right) \frac{\rho_i\rho_j}{\rho_\xi} - \frac{1}{\rho_\xi}\left(\rho^{l\bar k}\rho_l\rho_{j\bar k i} + \rho_{ij}\right).\\
\Gamma_{ij}^k & = \sum_{k=1}^n g^{k\bar l}\dd{z_i} g_{j\bar l} + g^{k\bar\xi}\dd{z_i} g_{j\bar\xi}\\
& = \frac{1}{f}\rho^{k\bar l}(f_\rho \rho_i \rho_{j\bar l} + f\rho_{j\bar l i} + f_{\rho\rho}\rho_i\rho_j\rho_{\bar l} + f_\rho\rho_{ij}\rho_{\bar l} + f_\rho \rho_j \rho_{i\bar l})\\
& \quad - \frac{1}{f\rho_{\bar\xi}} \rho^{k\bar l}\rho_{\bar l} (f_{\rho\rho} \rho_i \rho_j \rho_{\bar \xi} + f_\rho \rho_{ij} \rho_{\bar\xi})\\
& = \frac{f_\rho}{f} (\rho_i \delta_j^k + \rho_j \delta_i^k) + \rho^{k\bar l}\rho_{j\bar li}.
\end{align*}
\end{proof}
Let's state and prove our splitting lemma.
\begin{lemma}\label{splitting}
Let $M = \mathbb{P}(\mathcal{O}_\Sigma \oplus L)$ be the projective bundle such that the triple $(\Sigma, L, [\omega_0])$ satisfies the assumptions stated in P.\pageref{main}. Let $(M, \omega_t)$, $t \in [0,T)$ be the K\"ahler-Ricci flow $\partial_t \omega_t = -\Ric(\omega_t)$ with initial K\"ahler class $[\omega_0]$. Let $(x_i, t_i) \in M \times [0, T)$ be a sequence such that $t_i \to T$ and $K_i := \|\Rm(x_i)\|_{g(t_i)} \to \infty$ as $i \to \infty$. Define $g_i(t)$ to be rescaled dilated sequence by $K_i$ and $t_i$, i.e.
\[g_i(t) := K_i g(t_i + K^{-1}_i t), \quad t \in [-\beta_i, \alpha_i]\]
where $\beta_i \to \infty$, $\alpha_i \geq 0$ and $\alpha_i \to A \in [0, \infty]$. Suppose the curvature tensor of $g_i(t)$, $t \in [-\beta_i, \alpha_i]$ is uniformly bounded independent of $i$, i.e. there exists $C > 0$ independent of $i$ such that
\[\sup_{M \times [-\beta_i, \alpha_i]} \|\Rm\|_{g_i(t)} \leq C.\]
Then, after passing to a subsequence, $(M^{n+1}, g_i(t), x_i)$ converges smoothly in pointed Cheeger-Gromov sense to a complete ancient K\"ahler-Ricci flow $(M_\infty, g_\infty(t), x_\infty)$ whose universal cover is of the form
\[(N_1^n \times N_2^1, h_1(t) \oplus h_2(t)), \quad t \in (-\infty, A]\]
where $(N_i, h_i(t))$, $i = 1, 2$, are K\"ahler-Ricci flow solutions.
\end{lemma}
\begin{proof}
By the uniform boundedness condition of $\|\Rm\|_{g_i(t)}$ over $M \times [-\beta_i, \alpha_i]$, the subsequential Cheeger-Gromov convergence can be done by Hamilton's compactness theorem and Perelman's local non-collapsing theorem. See \cite{Chow1, Ham95,Perl1}, etc. Furthermore, we may assume the complex structure of $J$ of $M$ converges after passing to a subsequence to a complex structure $J_\infty$ of $M_\infty$. That makes $(M_\infty, J_\infty)$ K\"ahler because $\nabla^{g_\infty} J_\infty = \lim_{i\to\infty} \nabla^{g_i} J = 0$.

We will use the well-known de Rham's holonomy splitting theorem, which asserts that if the tangent bundle $TM_\infty$ admits an irreducible decomposition $\bigoplus_{i=1}^k E_i$ under the holonomy group action, i.e. parallel translation, then the universal cover of $M_\infty$ splits isometrically as $(M_\infty, g) = \prod_{i=1}^k N_i^{\dim E_i}$ with $TN_i^{\dim E_i} = E_i$. Note that in the K\"ahler case where the holonomy group is a subgroup of the unitary group, each $N_i$ is also K\"ahler.

Suppose $(M_\infty, g_\infty(t), x_\infty)$ is the pointed Cheeger-Gromov limit obtained above. We would like to show it (precisely, the universal cover) splits isometrically into a product. According to the nature of the collapsing of the $\mathbb{P}^1$-fiber, it is natural to guess that one factor of the split product should correspond to the base and the other should correspond to the fiber. Based on these, we define the following unit vector fields
\begin{align}
Z^j_{g_i(t)} & := \frac{1}{\|\dd{z_j}\|_{g_i(t)}} \dd{z_j} = \frac{1}{\sqrt{K_i (f\rho_{j\bar j}+f_\rho|\rho_j|^2)}} \dd{z_j},\\
\Xi_{g_i(t)} & := \frac{1}{\sqrt{K_i f_\rho}\rho_\xi} \dd{\xi}.
\end{align}
Then we have $\|Z^j_{g_i(t)}\|_{g_i(t)} = \|\Xi_{g_i(t)}\|_{g_i(t)} = 1$. After passing to a subsequence, they converge to vector fields $Z^j_{g_\infty(t)}$ and $\Xi_{g_\infty(t)}$ in the limit $M_{\infty}$. 

We will show that the real distribution $E_{\infty} = \text{span}_\R \{\Re(\Xi_{g_\infty(t)}), \Im(\Xi_{g_\infty(t)})\}$ is invariant under parallel translation. Here $\Re$ and $\Im$ denote denote the real and imaginary parts respectively. For simplicity, we will denote $Z^j_{g_i(t)}$ as $Z^j_i$ and $\Xi_{g_i(t)}$ as $\Xi_i$ for any $i \in \N \cup \{\infty\}$.

It is worthwhile to note that $E_{\infty}^{\perp} = \text{span}_\R \{\Re(Z^j_\infty), \Im(Z^j_\infty)\}_{j=1}^n$ since 
\begin{equation*}
\abs{\langle \Xi_i, \bar{Z}^j_i \rangle_{g_i(t)}} = \abs{\frac{1}{\sqrt{K_i f_\rho}\rho_\xi} \cdot \frac{1}{\sqrt{K_i g_{j\bar j}}} \cdot K_i f_\rho \rho_\xi \rho_{\bar j}} \leq \sqrt{\frac{f_\rho}{f}} \cdot \frac{\abs{\rho_j}}{\sqrt{\rho_{j\bar j}}},
\end{equation*}
which tends to $0$ as $i \to \infty$ using Lemma \ref{est_2}. Note that we have used $g_{j\bar j} = f\rho_{j\bar j} + f_\rho |\rho_j|^2 \geq f\rho_{j\bar j}$. Since $\rho = \log |\xi|^2 + \log \phi(\mathbf{z})$ by \eqref{rho_coord}, the term $\frac{\rho_j}{\sqrt{\rho_{j\bar j}}}$ is independent of $\xi$, $i$ and $t$, and hence is uniformly bounded near $\rho=\pm\infty$. Therefore $\Xi_\infty$ is orthogonal to each of $Z^j_\infty$, i.e. $E_{\infty}^{\perp} = \text{span}_\R \{\Re(Z^j_\infty), \Im(Z^j_\infty)\}_{j=1}^n$.

In order to show $E_\infty$ is invariant under parallel translation, we need to show that by parallel translating $\Xi_\infty$ along any vector field $X$ on $M_\infty$, it stays inside $E_\infty$. We will prove it by showing $\nabla^{\infty}_X \Xi_\infty$ lies inside $E_\infty$, or equivalently, orthogonal to $E_\infty^{\perp}$. We will make use of the Christoffel symbols calculated in Lemma \ref{christoffel}, 
\begin{align*}
\nabla_{\Xi_i}\Xi_i & = \frac{1}{\sqrt{K_i f_\rho} \rho_\xi} \nabla_{\xi} \left(\frac{1}{\sqrt{K_i f_\rho} \rho_\xi} \dd{\xi}\right)\\
& = \frac{1}{\sqrt{K_i f_\rho} \rho_\xi} \left(\dd{\xi}\left(\frac{1}{\sqrt{K_i f_\rho}\rho_\xi}\right)\dd{\xi} + \frac{1}{\sqrt{K_i f_\rho}\rho_\xi}\left(\Gamma_{\xi\xi}^\xi \dd{\xi} + \Gamma_{\xi\xi}^j \dd{z_j}\right)\right)\\
& = \frac{1}{\sqrt{K_i f_\rho} \rho_\xi} \left(\frac{1}{\sqrt{K_i}} \left(\frac{1}{\sqrt{f_\rho}} - \xi \frac{f_{\rho\rho}}{2f_\rho^{3/2}}\rho_\xi\right) \dd{\xi} + \frac{1}{\sqrt{K_i f_\rho}\rho_\xi} \left(\frac{f_{\rho\rho}}{f_\rho} - 1\right)\rho_\xi \frac{\partial}{\partial \xi}\right)\\
& = \frac{1}{K_i \sqrt{f_\rho}\rho_\xi} \left(\frac{1}{\sqrt{f_\rho}}\dd{\xi} - \frac{f_{\rho\rho}}{2\sqrt{f_\rho}f_\rho}\dd{\xi} + \frac{1}{\sqrt{f_\rho}} \frac{f_{\rho\rho}}{f_\rho}\dd{\xi} - \frac{1}{\sqrt{f_\rho}}\dd{\xi}\right)\\
& = \frac{\xi}{2K_i f_\rho} \left(\frac{f_{\rho\rho}}{f_\rho}\right) \dd{\xi}.
\end{align*}
Taking inner product with the vectors along the base direction, we have
\begin{align*}
\langle \nabla_{\Xi_i}\Xi_i, \bar{Z}^j_i \rangle_{g_i} & = \frac{\xi}{2K_i f_\rho} \cdot \frac{f_{\rho\rho}}{f_\rho} \cdot \frac{1}{\sqrt{K_i(f\rho_{j\bar j}+f_\rho|\rho_j|^2)}} \cdot K_i f_\rho \rho_\xi \rho_{\bar j},\\
|\langle \nabla_{\Xi_i}\Xi_i, \bar{Z}^j_i \rangle_{g_i}| & \leq \frac{1}{2\sqrt{K_i}} \cdot \frac{f_{\rho\rho}}{f_\rho} \cdot \frac{|\rho_{\bar j}|}{\sqrt{f \rho_{j\bar j}}}.
\end{align*}
Letting $i \to \infty$, we get $\langle \nabla_{\Xi_\infty}\Xi_\infty, \bar{Z}^j_\infty \rangle_{g_\infty} = 0$ for any $j = 1,\ldots, n$, here $\nabla$ is the Levi-Civita connection with respect to $g_\infty$. We have used the estimates proved in Lemma \ref{est_2}, which says $f_{\rho\rho}/f_\rho = O(1)$ and $f = O(1)$, as well as the fact that $K_i \to \infty$. This proves $\nabla_{\Xi_\infty} \Xi_\infty \in E_\infty$.

Similarly, by parallel translating $\Xi_\infty$ along $\bar{\Xi}_\infty$, we calculate
\begin{align*}
\nabla_{\bar{\Xi}_i} \Xi_i & = \frac{1}{\sqrt{K_i f_\rho}\rho_{\bar\xi}} \nabla_{\bar\xi} \left(\frac{1}{\sqrt{K_i f_\rho}} \xi \dd{\xi}\right)\\
& = -\frac{1}{2 K_i f_\rho} \left(\frac{f_{\rho\rho}}{f_\rho}\right) \xi \dd{\xi}\\
& = -\nabla_{\Xi_i}\Xi_i.
\end{align*}
Hence, we also have $|\langle \nabla_{\bar{\Xi}_i} \Xi_i, \bar{Z}^j_i \rangle_{g_i(t)}| \to 0$ 
for any $j = 1, \ldots, n$ as $i \to \infty$ and that proves $\langle \nabla_{\bar\Xi_\infty} \Xi_\infty \bar{Z}_\infty^j \rangle_{g_\infty} = 0$ and so $\nabla_{\bar\Xi_\infty} \Xi_\infty \in E_\infty$.

The other calculations are similar:
\begin{align*}
\nabla_{Z^j_i} \Xi_i & = \frac{1}{\sqrt{K_i g_{j\bar j}}} \nabla_j \left(\frac{1}{\sqrt{K_i f_\rho}\rho_\xi} \dd{\xi}\right)\\
& = \frac{1}{\sqrt{K_i g_{j\bar j}}} \left(\dd{z_j} \left(\frac{1}{\sqrt{K_i f_\rho} \rho_\xi}\right) \dd{\xi} + \frac{1}{\sqrt{K_i f_\rho} \rho_\xi} \left(\Gamma_{j\xi}^k \dd{z_k} + \Gamma_{j\xi}^{\xi} \dd{\xi}\right)\right)\\
& = \frac{1}{K_i \sqrt{g_{j\bar j}}} \left(-\frac{\rho_j}{2\sqrt{f_\rho}\rho_\xi}\left(\frac{f_{\rho\rho}}{f_\rho}\right) \dd{\xi}+\frac{1}{\sqrt{f_\rho}\rho_\xi} \left(\frac{f_\rho}{f}\delta_j^k \rho_\xi \dd{z_k} + \left(\frac{f_{\rho\rho}}{f_\rho}-\frac{f_\rho}{f}\right)\rho_j \dd{\xi}\right)\right)\\
& = \frac{1}{K_i \rho_\xi \sqrt{f_\rho g_{j\bar j}}} \left(\left(\frac{f_{\rho\rho}}{2f_\rho} - \frac{f_\rho}{f}\right) \rho_j \dd{\xi} + \frac{f_\rho}{f} \rho_\xi \dd{z_j}\right),
\end{align*}
\begin{align*}
\langle \nabla_{Z^j_i} \Xi_i, \bar{Z}^k_i \rangle_{g_i(t)} & = \frac{1}{K_i \rho_\xi \sqrt{f_\rho g_{j\bar j}}\sqrt{K_i g_{k\bar k}}}\\
& \quad \times \left\{\left(\frac{f_{\rho\rho}}{2f_\rho} - \frac{f_\rho}{f}\right) \rho_j K_i f_\rho \rho_\xi \rho_{\bar k} + \frac{f_\rho}{f} \rho_\xi K_i (f \rho_{j\bar k} + f_\rho \rho_j \rho_{\bar k}) \right\}\\
& = \sqrt{\frac{f_\rho}{K_i g_{j\bar j} g_{k\bar k}}} \left(\frac{f_{\rho\rho}}{2f_\rho} \rho_j \rho_{\bar k} + \rho_{j\bar k}\right).
\end{align*}
Hence
\begin{align*}
|\langle \nabla_{Z^j_i} \Xi_i, \bar{Z}^k_i \rangle_{g_i(t)}| & \leq \sqrt{\frac{f_\rho}{K_i}} \frac{1}{f} \left(\abs{\frac{f_{\rho\rho}}{2f_\rho} \cdot \frac{\rho_j \rho_{\bar k}}{\sqrt{\rho_{j\bar j}}\sqrt{\rho_{k\bar k}}}} + \abs{\frac{\rho_{j\bar k}}{\sqrt{\rho_{j\bar j}}\sqrt{\rho_{k\bar k}}}}\right)
\end{align*}
and so $|\langle \nabla_{Z^j_i} \Xi_i, \bar{Z}^k_i \rangle_{g_i(t)}| \to 0$ as $i \to \infty$ since by Lemma \ref{est_2} we have $f_\rho = O(T-t)$.

Finally, we have
\begin{align*}
\nabla_{\bar{Z}^j_i} \Xi_i & = \frac{1}{\sqrt{K_i g_{j\bar j}}} \dd{\bar{z}_j} \left(\frac{1}{\sqrt{K_i f_\rho}\rho_\xi} \dd{\xi} \right) \\
& = -\frac{1}{K_i f_\rho \sqrt{g_{j\bar j}}} \left(\frac{f_{\rho\rho}}{f_\rho}\right) \rho_{\bar j} \xi \dd{\xi},
\end{align*}
\begin{equation*}
\langle \nabla_{\bar{Z}^j_i} \Xi_i, \bar{Z}^k_i \rangle_{g_i(t)} = -\frac{1}{\sqrt{K_i g_{j\bar j} g_{k\bar k}}} \left(\frac{f_{\rho\rho}}{f_\rho}\right) \rho_{\bar j}\rho_{\bar k},
\end{equation*}
\begin{equation*}
|\langle \nabla_{\bar{Z}^j_i} \Xi_i, \bar{Z}^k_i \rangle_{g_i(t)}| \leq \frac{1}{f \sqrt{K_i}} \left(\frac{f_{\rho\rho}}{f_\rho}\right) \cdot \abs{\frac{\rho_j \rho_k}{\sqrt{\rho_{j\bar j} \rho_{k\bar k}}}}.
\end{equation*}
Hence $|\langle \nabla_{\bar{Z}^j_i} \Xi_i, \bar{Z}^k_i \rangle_{g_i(t)}| \to 0$ as $i \to \infty$.

Since $\{Z^j_{\infty}, \Xi_\infty\}_{j=1}^n$ spans the whole $T_\C M_\infty$, the above calculations show that for any vector field $X$ on $(M_{\infty}, g_{\infty}(t))$, one has $\langle \nabla_X \Xi_\infty, \bar{Z}^j_\infty \rangle_{g_\infty} = 0$ for any $j = 1, 2, \ldots, n$. Therefore, $\nabla_X \Re(\Xi_\infty), \nabla_X \Im(\Xi_\infty) \in E_\infty$. This shows whenever we have $V_x \in \eval{E_{\infty}}{x}$, $x \in M_\infty$ and let $V(s) \in TM$ be the parallel translation of $V_x$ along a curve $\gamma(s)$, then $V(s) \in E_\infty$. To see this, write $V(s) = V^T(s) + V^\perp(s)$ where $V^T(s) \in E_{\infty}$ and $V^\perp(s) \in E_{\infty}^{\perp}$ for any $s$. By the above calculation, we have $\nabla_{\gamma'(s)} V^T(s) \in E_{\infty}$ for any $s$. Therefore,
\[0 = \nabla_{\gamma'(s)} V(s) = \nabla_{\gamma'(s)} V^T(s) + \nabla_{\gamma'(s)} V^\perp(s).\]
Hence $\nabla_{\gamma'(s)} V^\perp(s)$ also lies inside $E_{\infty}$. By the fact that $V^\perp(s) \perp E_{\infty}$, we have
\begin{equation*}
\frac{d}{ds} \|V^\perp(s)\|^2 = 2\left\langle \nabla_{\gamma'(s)} V^\perp(s), V^\perp(s) \right\rangle = 0.
\end{equation*}
It implies that $\|V^\perp(s)\| \equiv \|V^\perp(0)\| = 0$ for any $s$. In other words, $V(s) \equiv V^T(s) \in E_{\infty}$ for any $s$. Therefore, $E_{\infty}$ is invariant under parallel transport. By the de Rham's decomposition theorem, our splitting lemma follows.
\end{proof}

\section{Singularity Analysis}
The splitting lemma in the previous section allows a dimension reduction for our singularity analysis. The ultimate goal of this section is to analyze the singularity formation of the Ricci flow on our projective bundles $M = \mathbb{P}(\mathcal{O}_\Sigma \oplus L)$ whose $\mathbb{P}^1$-fiber collapses near the singularity. We are going to prove that the K\"ahler-Ricci flow $(M, g(t))$ must be of Type I (see definition below) and the singularity model is $\C^n \times \mathbb{P}^1$, in a sense that one can choose a sequence $(x_i, t_i)$ in space-time in the high curvature region such that the universal cover of the Cheeger-Gromov limit of the rescaled dilated sequence is isometric to $(\C^n \times \mathbb{P}^1, \|d\mathbf{z}\|^2 \oplus \omega_{\text{FS}}(t))$.
Here $\omega_{\text{FS}}(t)$ is the shrinking Fubini-Study metric.

According to the blow-up rate of the Riemann curvature tensor, the singularity type of a Ricci flow solution which encounters finite-time singularity is classified as in \cite{Ham95}.
\begin{definition}
Let $(M, g(t))$ be a Ricci flow solution $\partial_t g(t) = -\Ric(g(t))$ on a closed manifold $M$ which becomes singular at a finite time $T$. We call the Ricci flow encounters
\begin{itemize}
\item \textbf{Type I singularity} if $\sup_{M \times [0,T)}(T-t)\|\Rm\|_{g(t)} < \infty$;
\item \textbf{Type II singularity} if $\sup_{M \times [0,T)}(T-t)\|\Rm\|_{g(t)} = \infty$.
\end{itemize}
\end{definition}

We would like to remark that although the Type I/II classification of finite-time singularity was proposed in the early 90's, surprisingly the first compact Type II solution was constructed by Gu-Zhu in \cite{GuZh07} only recently in 2007.

In order to understand the singularity formation, we need to bring curvatures into the topic. Therefore, we will compute and analyze the Riemann curvature tensor of our projective bundle $M$ which is equipped with momentum profile $f$. Recall that for K\"ahler manifolds, the Riemann curvature $(3,1)$-tensor can be computed using the formula
\begin{equation*}
R_{A\bar BC}^D = -\dd{\bar{z}^B} \Gamma_{AC}^D
\end{equation*}
where $A, B, C, D = 1, \ldots, n$ or $\xi$. The non-zero components of the Riemann curvature tensor are given below. For the ease of inspection of the norm $\|\Rm\|$ later on, we will split the components into five groups according to the number of $\xi$-indexes.
\label{riem_comp}
\begin{align*}
R_{i\bar jk}^l & = -(\log f)_{\rho\rho} \rho_{\bar j} (\rho_i \delta_{kl} + \rho_k \delta_{il}) - (\log f)_\rho (\delta_{ij}\delta_{kl} + \delta_{jk}\delta_{il}) - (\rho^{l\bar p}\rho_{ik\bar p})_{\bar j}. \\ \\
R_{i\bar\xi k}^ l & = -(\log f)_{\rho\rho} \rho_{\bar\xi} (\rho_i \delta_{kl} + \rho_k \delta_{il}),\\
R_{i\bar j\xi}^l & = -(\log f)_{\rho\rho}\rho_{\bar j}\rho_\xi \delta_{il},\\
R_{i\bar jk}^\xi & = -\frac{1}{\rho_\xi} (\log f_\rho - 2\log f)_{\rho\rho}\rho_{\bar j}\rho_i\rho_k, \\
& \quad - \frac{1}{\rho_\xi} (\log f_\rho - 2\log f)_\rho (\rho_{i\bar j} \rho_k + \rho_{k\bar j} \rho_i) + \frac{1}{\rho_\xi} (\rho^{l\bar p} \rho_l \rho_{i\bar p k} + \rho_{ik})_{\bar j},\\
R_{\xi\bar jk}^l & = -(\log f)_{\rho\rho} \rho_{\bar j}\rho_\xi\delta_{kl}.\\ \\
R_{i\bar\xi k}^\xi & = -(\log f_\rho - 2\log f)_{\rho\rho} \frac{\rho_{\bar\xi}}{\rho_\xi}\rho_i\rho_k,\\
R_{i\bar\xi\xi}^l & = -(\log f)_{\rho\rho}|\rho_\xi|^2 \delta_{ik},\\
R_{\xi\bar\xi k}^l & = -(\log f)_{\rho\rho}|\rho_\xi|^2\delta_{kl},\\
R_{\xi\bar j\xi}^l & = 0,\\
R_{\xi\bar jk}^\xi & = -(\log f_\rho - \log f)_{\rho\rho}\rho_{\bar j}\rho_k - (\log f_\rho - \log f)_\rho \rho_{k\bar j},\\
R_{l\bar j\xi}^\xi & = -(\log f_\rho - \log f)_{\rho\rho} \rho_{\bar j} \rho_l - (\log f_\rho - \log f)_\rho \rho_{l\bar j}.\\ \\
R_{l\bar\xi\xi}^\xi & = -(\log f_\rho - \log f)_{\rho\rho} \rho_{\bar\xi} \rho_i,\\
R_{\xi\bar\xi\xi}^l & = 0,\\
R_{\xi\bar\xi k}^\xi & = -(\log f_\rho - \log f)_{\rho\rho}\rho_{\bar\xi}\rho_k,\\
R_{\xi\bar j\xi}^\xi & = -(\log f_\rho)_{\rho\rho}\rho_{\bar j}\rho_\xi.\\ \\
R_{\xi\bar\xi\xi}^\xi & = -(\log f_\rho)_{\rho\rho}|\rho_\xi|^2.
\end{align*}
Since the understanding of $\|\Rm\|$ is crucial in analyzing the singularity according their type (I or II), we need an organized expression of $\|\Rm\|$ that is written in terms of our momentum profile $f$. Obviously, it would take loads of unnecessary work. However, in order to study the singularity model in our class of manifolds, it suffices to understand the asymptotics of $\|\Rm\|^2$ in terms of $f$ and its derivatives.
Recall from Lemma \ref{est_2} that $f = O(1)$, $\frac{1}{f} = O(1)$, $f_{\rho\rho}/f_\rho = O(1)$. Therefore we have the following asymptotics
\begin{align*}
(\log f)_\rho & = \frac{f_\rho}{f} = O(f_\rho),\\
(\log f)_{\rho\rho} & = \frac{f_{\rho\rho}}{f} - \frac{f^2_\rho}{f^2} = O(f_\rho),\\
(\log f_\rho)_\rho & = \frac{f_{\rho\rho}}{f_\rho} = O(1).
\end{align*}
The asymptotic of $(\log f_\rho)_{\rho\rho}$ is not yet known because it involves the third $\rho$-derivative of $f$ which we have not derived.

Also, the local expressions of $g$ and $g^{-1}$ have the following asymptotics
\begin{align*}
g_{i\bar j} & = O(1), \\
g_{i\bar\xi} & = g_{\bar i\xi} = g_{\xi\bar\xi} = O(f_\rho),\\
g^{i\bar j} & = g^{i\bar\xi} = g^{\bar i\xi} = O(1),\\
g^{\xi\bar\xi} & = O(f_\rho^{-1}).
\end{align*}

We claim that the norm $\|\Rm\|^2$ can be expressed in the following asymptotic form
\begin{lemma}
\begin{align}\label{riem}
\|\Rm\|^2_{g(t)} & = f^{-2}_\rho (\log f_\rho)^2_{\rho\rho} + O(f^{-1}_\rho (\log f_\rho)_{\rho\rho}) \\ 
& \quad + O(f^{-1}_\rho (\log f_\rho)^2_{\rho\rho}) + O((\log f_\rho)^2_{\rho\rho}) \notag\\
& \quad + O((\log f_\rho)_{\rho\rho}) + O(1). \notag
\end{align}
\end{lemma}
\begin{proof}
A generic term in $\|\Rm\|^2$ can be expressed as
\[g_{A\bar B}g^{C\bar D}g^{E\bar F}g^{G\bar H}R^A_{C\bar{F}G}\barrr{R^B_{D\bar{E}H}} \tag{**} \label{term}\] where $A, \ldots, H \in \{1, \ldots, n, \xi\}$. From Lemma \eqref{est_2}, we know $f_\rho = O(T-t)$ and so $f^{-1}_\rho$ is a bad term as it diverges as $t \to T$. The only factor in \eqref{term} which can contribute to a $f_\rho^{-1}$ is $g^{\xi\bar\xi}$, and there are at most three $g^{\xi\bar\xi}$'s in (\ref{term}). We are going to check that
\begin{enumerate}[(1)]
\item whenever $f^{-1}_\rho$ appears in \eqref{term} exactly once, there must at least one factor of $(\log f_\rho)_{\rho\rho}$ from the curvature components;
\item whenever $f^{-2}_\rho$ appears in \eqref{term}, there must be a $(\log f_\rho)^2_{\rho\rho}$ factor from the curvature components;
\item it is impossible for $f^{-3}_\rho$ to appear in \eqref{term}.
\end{enumerate}
Combining these, it is not difficult to see $\|\Rm\|^2$ satisfies the asymptotic form \eqref{riem}.

We start by arguing (1). Suppose there is exactly one $f^{-1}_\rho$ factor in \eqref{term}, we can assume WLOG that either $(C, D) = (\xi, \xi)$ or $(E, F) = (\xi, \xi)$. Suppose the former, we can check from the table of Riemann curvatures in P.\pageref{riem_comp} that almost all $R_{\xi\bar FG}^A$ terms have either asymptotics $O(f_\rho)$ (which cancels out $f^{-1}_\rho$) or a $(\log f_\rho)_{\rho\rho}$ factor. There is only one exception: $R_{\xi\bar jk}^\xi$ which has an $O(1)$-term from $(\log f_\rho)_\rho$. However, if both of $R^A_{C\bar{F}G}$ and $R^B_{D\bar{E}H}$ are taken to be in this form, then \eqref{term} becomes
\[g_{\xi\bar\xi} g^{\xi\bar\xi} g^{p\bar j} g^{k\bar q} R_{\xi\bar jk}^\xi \barrr{R_{\xi\bar pq}^\xi},\]
where the $g_{\xi\bar\xi} = O(f_\rho)$ cancels out the undesirable $f^{-1}_\rho$ factor, and end up with no $f_\rho^{-1}$ at all. Similar argument applies to the case $(E, F) = (\xi, \xi)$, and (1) is proved.

For (2), since $g^{\xi\bar\xi}$ is the only possible contribution to $f^{-1}_\rho$, at least two of $C, F, G$ (and their corresponding two of $D, E, G$) must be $\xi$. Check again the table of Riemann curvature components in P. \pageref{riem_comp}, we see all the terms with two lower $\xi$-indexes must either of $O(f_\rho)$-type or has a $(\log f_\rho)_{\rho\rho}$ factor. It proves (2).

For (3), the only possible case for $f^{-3}_\rho$ to appear is that all of $(C, D)$, $(E, F)$ and $(G, H)$ are $(\xi, \xi)$. The only possible choice for the curvature components are $R_{\xi\bar\xi\xi}^l$ and $R_{\xi\bar\xi\xi}^\xi$. However, the former is 0. For the latter case, all indexes will be $\xi$ and \eqref{term} becomes
\[g_{\xi\bar\xi}g^{\xi\bar\xi}g^{\xi\bar\xi}g^{\xi\bar\xi}R_{\xi\bar\xi\xi}^\xi\barrr{R_{\xi\bar\xi\xi}^\xi}\]
which can be computed easily as $f_{\rho}^{-2} (\log f_\rho)^2_{\rho\rho}$.

Finally, we remark that $g_{\xi\bar\xi}g^{\xi\bar\xi}g^{\xi\bar\xi}g^{\xi\bar\xi}R_{\xi\bar\xi\xi}^\xi\barrr{R_{\xi\bar\xi\xi}^\xi}$ is the only term that $f^{-2}_\rho (\log f_\rho)^2_{\rho\rho}$ appears, thanks to the fact that $R_{\xi\bar\xi\xi}^i = 0$. As a result, the leading term of \eqref{riem} is $f^{-2}_\rho (\log f_\rho)^2_{\rho\rho}$ with coefficient $1$ which can be easily verified by computing $g_{\xi\bar\xi}g^{\xi\bar\xi}g^{\xi\bar\xi}g^{\xi\bar\xi}R_{\xi\bar\xi\xi}^\xi\barrr{R_{\xi\bar\xi\xi}^\xi}$.
\end{proof}
Having understood the asymptotics of $\|\Rm\|^2$, we are in a position to study the singularity models. Let's first consider the Type I case:
\begin{theorem}\label{type_I}
Let $M = \mathbb{P}(\mathcal{O}_\Sigma \oplus L)$ be the projective bundle with the triple $(\Sigma, L, [\omega_0])$ satisfying the conditions listed in P.\pageref{main}. Let $(M, \omega_t)$ be the K\"ahler-Ricci flow $\partial_t \omega_t = -\Ric(\omega_t)$, $t \in [0,T)$ with initial K\"ahler class $[\omega_0]$. Suppose the flow encounters Type I singularity, then choose $(x_i, t_i)$ in space-time such that $K_i := \|\Rm(x_i, t_i)\|_{g(t_i)} = \max_M \|\Rm\|_{g(t_i)}$ and $t_i \to T$. Consider the rescaled dilated sequence of metrics $g_i(t) := K_i g(t_i + K^{-1}_i t)$, $t \in [-t_i K_i, (T-t_i)K_i)$. Then the pointed sequence $(M, g_i(t), x_i)$ converges, after passing to a subsequence, smoothly in pointed Cheeger-Gromov sense to an ancient $\kappa$-solution $(M_\infty, g_\infty(t), x_\infty)$, whose universal cover splits isometrically as
\[(\C^n \times \mathbb{P}^1, \|d\mathbf{z}\|^2 \oplus \omega_{\text{FS}}(t) ),\]
where $\|d\mathbf{z}\|^2$ is the Euclidean metric and $\omega_{\text{FS}}(t)$ denotes the shrinking Fubini-Study metric.
\end{theorem}
\begin{proof}
Suppose $C = C(n)$ is a constant depending only on $n$ such that $|R| \leq C(n)\|\Rm\|$. Since the blow-up factor $K_i$ is defined by $K_i = \max_M \|\Rm\|_{g(t_i)} = \|\Rm(x_i)\|_{g(t_i)}$, the scalar curvature at time $t_i$ satisfies $|R(g(t_i))| \leq CK_i$ on $M$.

One can compute the scalar curvature explicitly:
\begin{align*}
R_{g(t)} & = \Tr_{\omega_t} \Ric(\omega_t)\\
& = \frac{n(\nu - F_\rho)}{f} - \frac{F_{\rho\rho}}{f_\rho},
\end{align*}
where $F = \log f_\rho + n \log f$. Hence, 
\begin{equation*}
R_{g(t)} = - \frac{1}{f_\rho} (\log f_\rho)_{\rho\rho} + O(1).
\end{equation*}
Therefore, for any $\rho\in [-\infty, \infty]$ at $t_i$, we have
\begin{align*}
\abs{-\frac{1}{f_\rho}(\log f_\rho)_{\rho\rho} + O(1)} & \leq CK_i,\\
\abs{\frac{-1}{K_i f_\rho} (\log f_\rho)_{\rho\rho} + O(K^{-1}_i)} & \leq C.
\end{align*}
Recall that $K_i \to \infty$. Letting $i \to \infty$ yields
\begin{equation}\label{limit_ine}
\limsup_{i\to\infty} \abs{K_i^{-1} f^{-1}_\rho (\log f_\rho)_{\rho\rho}}_{(\rho, t_i)} \leq C.
\end{equation}

By considering the asymptotic expression of $\|\Rm\|^2$ given by \eqref{riem}, we have for any $\rho \in [-\infty, \infty]$ at time $t_i$,
\begin{align*}
1 \geq K^{-2}_i \|\Rm\|^2_{g(t_i)} & = (K_i f_\rho)^{-2} (\log f_\rho)^2_{\rho\rho} + O(K^{-2}_i f^{-1}_\rho (\log f_\rho)^2_{\rho\rho}) \\
& \quad + O(K_i^{-2}f_\rho^{-1} (\log f_\rho)_{\rho\rho}) + O(K_i^{-2} (\log f_\rho)^2_{\rho\rho})\\
& \quad + O(K_i^{-2} (\log f_\rho)_{\rho\rho}) + O(K_i^{-2}),
\end{align*}
where equality is achieved at $x_i$.

Letting $i \to \infty$ and using \eqref{limit_ine} and the fact that $f_\rho = O(T-t)$ from Lemma \ref{est_2}, we can deduce:
\begin{align}
\limsup_{i\to\infty} (K_i f_\rho)^{-2} (\log f_\rho)^2_{\rho\rho} & \leq 1, \quad \rho \in [-\infty, \infty], \quad t = t_i,\notag\\
\lim_{i\to\infty} \eval{(K_i f_\rho)^{-2} (\log f_\rho)^2_{\rho\rho}}{(x_i, t_i)} & = 1. \label{limit_eq}
\end{align}

Recall that $g_i(t) = K_i g(t_i + K^{-1}_i t)$, we then have
\begin{equation*}
R_{g_i(t)} = \eval{-\frac{1}{K_i f_\rho} (\log f_\rho)_{\rho\rho} + O(K_i^{-1})}{t_i + K^{-1}_i t}.
\end{equation*}
Letting $i \to \infty$, we have
\begin{equation}\label{limit_scalar}
R_{g_{\infty}(t)} = -\lim_{i \to \infty} \eval{\frac{1}{K_i f_\rho} (\log f_\rho)_{\rho\rho}}{t_i + K^{-1}_i t}.
\end{equation}
By strong maximum principle, the scalar curvature of every ancient solution must be either identically zero or everywhere positive. In our case, \eqref{limit_eq} and \eqref{limit_scalar} together implies $R_{g_{\infty}(0)} = 1$ and hence $R_{g_\infty(t)} > 0$ on $M \times (-\infty, 0]$. By our splitting lemma \ref{splitting}, we know that the limit manifold $M_{\infty}$ splits isometrically as a product $N_1^n \times N_2^1$, such that $TN_1^n = \text{span}_\R \{\Re(Z_\infty^j), \Im(Z_\infty^j)\}_{j=1}^n$ and $TN_2^1 = \text{span}_\R \{\Re(\Xi_\infty), \Im(\Xi_\infty)\}$. As a result, the curvature tensors also split as $\Rm_{M_\infty} = \Rm_{N_1^n} \oplus \Ric_{N_2^1}$. Next, we would like to compute the curvatures of each factor. Again, for simplicity we denote $Z_{g_i(t)}^j$ by $Z_i^j$ and $\Xi_{g_i(t)}$ by $\Xi_i$
\begin{align}\label{riem_0}
& |\Rm_{g_i(t)} (Z_i^j, \bar{Z}_i^k, Z_i^l, \bar{Z}_i^p)|\\
& = |K_i \langle \Rm(Z_i^j, \bar{Z}_i^k) Z_i^l, \bar{Z}_i^p \rangle_{g(t_i+K^{-1}_i t)}| \notag\\
& \leq \abs{K_i \left(\frac{1}{\sqrt{K_i}}\right)^4 \frac{1}{\sqrt{f\rho_{j\bar j}}} \frac{1}{\sqrt{f\rho_{k\bar k}}} \frac{1}{\sqrt{f\rho_{l\bar l}}} \frac{1}{\sqrt{f\rho_{p\bar p}}} R_{j\bar kl\bar p}}\notag\\
& = \frac{1}{K_i} O(1) \to 0 \quad \text{ as } i \to \infty. \notag
\end{align}
Hence $\Rm_{N^n} = 0$. Similarly, we have
\begin{align}\label{ric_pos}
\Ric_{g_i(t)} (\Xi_i, \bar{\Xi}_i) & = \frac{1}{\sqrt{K_i f_\rho}} \frac{1}{\sqrt{K_i f_\rho}} \frac{1}{|\rho_\xi|^2} (-(n\log f + \log f_\rho)_{\rho\rho} |\rho_\xi|^2)\\
& = -\frac{1}{K_i f_\rho} (\log f_\rho)_{\rho\rho} + O(K_i^{-1}).\notag
\end{align}
By \eqref{limit_scalar} and positivity of $R_{g_\infty(t)}$, we know that $\Ric_{g_\infty(t)}(\Xi_\infty, \bar{\Xi}_\infty) > 0$.

Since the K\"ahler-Ricci flow $g(t)$ is of Type I, the ancient solution obtained by the blow-up sequence is also of Type I, i.e. $\sup_{M \times (-\infty, 0]} |t|\|\Rm\|_{g_\infty(t)} < \infty$, and is $\kappa$-non-collapsed. The limit solution splits as a product $(N_1^n, h_1(t)) \times (N_2^1, h_2(t))$ which we know $N_1^n$ is flat and $N_2^1$ has positive curvature. According to Hamilton's classification of ancient $\kappa$-solution \cite{Ham95} (see also \cite{CNL}), $(N_2^1, h_2(t))$ must be the shrinking round 2-sphere.

To conclude, if the K\"ahler-Ricci flow $(M, g(t))$ is of Type I, then the universal cover of the limit solution $(M_\infty, g_\infty(t))$ of the rescaled dilated sequence $g_i(t)$ is isometric to 
\[(\C^n \times \mathbb{P}^1, \|d\mathbf{z}\|^2 \oplus \omega_{\text{FS}}(t) ).\]

\end{proof}
Next, we will rule out the possibility of Type II singularity on $(M, g(t))$. We will show that by a standard point-picking argument for Type II singularity, one can form a rescaled dilated sequence of metrics which converges, after passing to a subsequence, to a product of the cigar soliton and a flat factor. By Perelman's local non-collapsing result, such limit model is not possible. Let's state this result and give its proof.
\begin{theorem}\label{no_type_II}
Let $M = \mathbb{P}(\mathcal{O}_\Sigma \oplus L)$ be the projective bundle with the triple $(\Sigma, L, [\omega_0])$ satisfying the conditions listed in P.\pageref{main}. Let $(M, \omega_t)$ be the K\"ahler-Ricci flow $\partial_t \omega_t = -\Ric(\omega_t)$, $t \in [0,T)$ with initial K\"ahler class $[\omega_0]$. Then $(M, g(t))$ must be of Type I, i.e. Type II singularity is not possible.
\end{theorem}
\begin{proof}
First take an increasing sequence $T_i \to T$. Let $(x_i, t_i) \in M \times [0, T_i]$ be such that
\begin{align*}
(T_i - t_i) \|\Rm\|(x_i, t_i) & = \max_{M\times[0,T_i]} (T_i - t)\|\Rm\|_{g(t)} \\
& = \max_{M \times [-K_i t_i, K_i (T_i - t_i)]} (T_i - (t_i + K^{-1}_i t)) \|\Rm\|_{g(t_i + K^{-1}_i t)}.
\end{align*}
We denote $K_i = \|\Rm\|(\rho_i, t_i)$, then $K_i (T_i - t_i) \to \infty$ by the Type II condition. 

As in the Type I case, we let $C=C(n)$ be a constant depending only $n$ such that $|R_{g(t)}| \leq C\|\Rm\|_{g(t)}$. Recall that scalar curvature has the following asymptotic expression:
\[R_{g(t)} = -\frac{1}{f_\rho}(\log f_\rho)_{\rho\rho} + O(1).\]
Hence for any $\rho \in [-\infty,\infty]$, $t \in [0, T_i]$, we have
\begin{align*}
\abs{-\frac{1}{f_\rho}(\log f_\rho)_{\rho\rho} + O(1)} & \leq \frac{C(T_i-t_i)K_i}{T_i-(t_i+K^{-1}_it)},\\
\abs{\frac{-1}{K_i f_\rho} (\log f_\rho)_{\rho\rho} + O(K^{-1}_i)} & \leq \frac{C(T_i-t_i)}{T_i-(t_i+K^{-1}_it)},
\end{align*}
where we evaluate the left-hand side at $t_i + K^{-1}_i t$. Letting $i\to\infty$, and using the fact that $\frac{(T_i - t_i) - K^{-1}_i t}{T_i - t_i} = 1 - \frac{t}{K_i (T_i - t_i)} \to 1$,
one can show
\begin{equation}\label{bound}
\limsup_{i \to \infty} |\eval{(K_i f_\rho)^{-1} (\log f_\rho)_{\rho\rho}}{(x, t_i +K^{-1}_i t)} \leq 1 \quad \text{ for any } (x,t).\\
\end{equation}
At $(x_i, t_i)$ we have $(T_i - t_i)^2 \|\Rm\|^2(x_i, t_i) = (T_i - t_i)^2 K_i^2$. Consider the asymptotic expression of $\|\Rm\|^2$ as in the Type I case, one can then show
\begin{equation}\label{equality}
\lim_{i\to\infty} \eval{\frac{1}{K_i^2 f^2_{\rho}} (\log f_\rho)^2_{\rho\rho}}{(x_i, t_i)} = 1.
\end{equation}

As $K_i \to \infty$, our splitting lemma \ref{splitting} also implies the limit solution $(M_\infty, g_\infty(t))$ splits isometrically as a product $(N_1^n \times N_2^1, h_1(t) \oplus h_2(t))$. As in the Type I case, $\Rm_{N_1^n}$ and $\Ric_{N_2^1}$ can be found by \eqref{riem_0} and \eqref{ric_pos}:
\begin{align*}
\Rm_{g_i(t)} (Z_i^j, \bar{Z}_i^k, Z_i^l, \bar{Z}_i^p) & = \frac{1}{K_i} O(1),\\
\Ric_{g_i(t)} (\Xi_i, \bar{\Xi}_i) & = -\frac{1}{K_i f_\rho} (\log f_\rho)_{\rho\rho} + O(K_i^{-1}).
\end{align*}
Letting $i \to \infty$, we have $\Rm_{N_1^n}(h_1(t)) = 0$ and
\begin{align}\label{max_ric_ach}
1 \geq \Ric_{N_2^1}(h_2(t)) & > 0 & \text{from \eqref{bound}},\\
\Ric_{N_2^1}(x_\infty, h_2(0)) & = h_2(0) & \text{from \eqref{equality}}.\notag
\end{align}

$(M_\infty, g_\infty(t))$ is an eternal solution to the K\"ahler-Ricci flow since we have $(T_i - t_i)K_i \to \infty$. By our splitting lemma, so does $(N_2^1, h_2(t))$. From \eqref{max_ric_ach}, the space-time maximum of the scalar curvature of $(N_2^1, h_2(t))$ is achieved at $(x_\infty, 0)$. Hence by Hamilton's classification of eternal solutions (see the Main Theorem of \cite{Ham93}), $(N_2^1, h_2(t))$ is a steady gradient soliton. In case of $\dim_\R = 2$, it must be the cigar soliton (see Section 26.3 of \cite{Ham95}). However, by Perelman's local non-collapsing \cite{Perl1}, the Cheeger-Gromov limit $(M_\infty, g(t))$ must be $\kappa$-non-collapsed at all scales, and so the product of cigar soliton and a flat space is not a possible singularity model. It leads to a contradiction and hence completes our proof.
\end{proof}

\begin{remark}
Throughout this paper we have focused on Case 1 and Case 2(i) in P.\pageref{cases}. We would like to point out as a final remark that for Case 2(iii) we expect one could mimic Section 5.2 in \cite{SgWk09} and also \cite{SgWk10, SgWk11} to show the contraction of $\Sigma_0$ near the singular time. For singularity models obtained by rescaling analysis in Case 2(iii), it is conjectured in \cite{FIK03} that for $(\Sigma, \omega_\Sigma) = (\mathbb{P}^n, \omega_{\FS})$ the singularity should be modelled on K\"ahler-Ricci solitons on $\mathcal{O}(-k)$-bundles over $\mathbb{P}^n$.
\end{remark}

\newcommand{\etalchar}[1]{$^{#1}$}
\providecommand{\bysame}{\leavevmode\hbox to3em{\hrulefill}\thinspace}
\providecommand{\MR}{\relax\ifhmode\unskip\space\fi MR }
\providecommand{\MRhref}[2]{%
  \href{http://www.ams.org/mathscinet-getitem?mr=#1}{#2}
}
\providecommand{\href}[2]{#2}

\end{document}